\documentclass{article}
\usepackage[utf8]{inputenc}
\usepackage[english]{babel}
\usepackage{authblk}
\usepackage{amsmath,amsthm,amssymb,graphicx,stmaryrd,enumerate,bbm,alltt}
\usepackage[bookmarksopen, bookmarksnumbered]{hyperref}
\newcommand{\nin}{\not\in}
\newcommand{\nobracket}{}
\newcommand{\nocomma}{}

\newcommand{\tmem}[1]{{\em #1\/}}

\newcommand{\tmop}[1]{\ensuremath{\operatorname{#1}}}

\newcommand{\tmverbatim}[1]{{\ttfamily{#1}}}
\newcommand{\R}{\mathbbm{R}}
\newcommand{\Z}{\mathbbm{Z}}
\newcommand{\N}{\mathbbm{N}}

\newtheorem{theorem}{Theorem}[section]
\newtheorem{lemma}[theorem]{Lemma}

\newtheorem{proposition}[theorem]{Proposition}
\newtheorem{corollary}[theorem]{Corollary}

\newtheorem{definition}[theorem]{Definition}
\newtheorem{remark}[theorem]{Remark}

\newenvironment{enumerateroman}{\begin{enumerate}[i.] }{\end{enumerate}}

\begin{document}
\markboth{Linjun Li}
{Shrinking braids and Left distributive monoid}

\title{Shrinking braids and Left distributive monoid}

\author{Linjun Li 
\thanks{
Department of Mathematics, University of Pennsylvania, Philadelphia, PA. e-mail: linjun@sas.upenn.edu}}
\date{}
\maketitle

\begin{abstract}
  We consider a natural generalization of braids which we call \emph{shrinking braids}. We state the relations of shrinking braids and use them to define algebraically the monoid $R$.
  We endow a subset of $R$ with a \emph{left distributive monoid}
  structure and use it to extend the Dehornoy order on $B_{\infty}$ to an order on
  $R$. By using this order, we prove that $R$ is isomorphic to the monoid which is generated (geometrically) by shrinking braids.
\end{abstract}

\section{Introduction}

{\tmem{Left distributive system (LD system) }}is a set $X$ with a
composition $\cdot$ on it such that the left distributive law holds. i.e. $a
\cdot ( b \cdot c ) = ( a \cdot b ) \cdot ( a \cdot c )$. LD systems emerge in
the study of knot invariance as fundamental quandles and show many interesting
properties. If there exists another composition $\circ$ on $X$ such that some
algebra laws (defined in Section \ref{sec:env}) between $\circ$ and $\cdot$ hold, we call
$X$ a \emph{left distributive monoid} (or \emph{LD monoid}).

The \emph{free monogenerated left distributive monoid} was naturally found
in the field of large cardinal in set theory (\cite{laver1992left}) as the algebra of elementary embeddings. There,
the composition $\cdot$ appears as the iteration of two elementary embeddings
and $\circ$ appears as the composition of elementary embeddings as functions.

Surprisingly, Dehornoy (\cite{dehornoy1994braid}) proved that there is a left distributive
composition $\cdot$ on $B_{\infty}$ (the infinite braid group) such that the
$\tmop{LD}$ system generated by the identity is free. This
discovery opened a project to seek relations between LD structure and low
dimensional topology. For example, this leads to the discovery of the Dehornoy
order on $B_{\infty}$ which shares interesting properties such as the
well-orderness on positive braids.

On the other hand, there are also many problems along this road since it is
still unclear to us the relation between these two representations of the free
LD system, namely, as elementary embeddings and as braids. A natural question
is whether one can construct LD monoid structures from the braid diagram or its natural extensions.

Along this question, in \cite{dehornoy1998transfinite} and \cite{dehornoy2006group}, two extensions of braid diagrams are constructed and extended LD structures are found from these diagrams. In \cite{dehornoy1998transfinite}, the \emph{transfinite braids} are constructed and the \emph{extended braid monoid} $EB_{\infty}$ is defined using the transfinite braids. It is proved that there is an $LD$ monoid structure on $EB_{\infty}$. In \cite{dehornoy2006group}, the \emph{parenthesized braid group} $B_{\bullet}$ is constructed as a ``mixture'' of the braid group and the Thompson group, and an augmented LD system (which is an
intermediate object between LD system and LD monoid) is found on it. 

Our main goal here is to study another extension of braids called \emph{shrinking braids} (see Section \ref{sec:shrink}). Shrinking braids belong to the set of two dimensional cobordisms whose input and output boundaries are two sequences of circles. Besides the braid-like cobordisms (Figure \ref{fig:braid}), the shrinking braids contain additionally pant-like cobordisms (Figure \ref{fig:pant}) which make the shrinking braids not invertible. In this way, shrinking braids are different from 
the parenthesized braids which are invertible. Meanwhile, shrinking braids do not contain transfinitely many input or output circles which are used in the construction of the transfinite braids.
On the other hand, the monoid $R$ (Definition \ref{def:define-R}) and LD monoid $\mathcal{B}$ (equation \eqref{eq:def_matc-B}) constructed from shrinking braids have strong connection with $B_{\bullet}$ and $EB_{\infty}$. See Remark \ref{rem:relation-to-parenthesized} and Remark \ref{rem:rel-EB} below.

The algebraic relations of shrinking braids are given in Propostition \ref{prop:relationofS} and these relations are proved to completely characterize shrinking braids (Theorem \ref{thm:main}).
It should be noticed that in the construction of shrinking braids, the cobordisms we considered are embedded into $\R^{3}$ and the equivalence between two cobordisms are given by certain ambient isotopy (Definition \ref{def:cobordisms}). This is similar to the case of braid diagrams and different from the usual definition of cobordisms where there is no ambient space. 
For results on two dimensional cobordisms with no ambient space and their relationship with algebra, see e.g. \cite{kock2004frobenius}.

\subsection*{Organization of remaining text}
In Section \ref{sec:shrink}, we use $2$-cobordisms to construct the geometric subject shrinking braids and state the relations of shrinking braids (Proposition \ref{prop:relationofS}). 
The next five sections are devoted to prove Theorem \ref{thm:main} which says that the monoid of shrinking braids are completely determined by the relations in Proposition \ref{prop:relationofS}. 
In order to do so, we define algebraically the monoid $R$ by these relations and study its basic properties (Section \ref{sec:basic}).
Then we identify the free monogenerated $LD$ monoid as the \emph{enveloping $LD$ monoid} of the free monogerated $LD$ system (Section \ref{sec:env}).
After that, in Section \ref{sec:monoidR}, we find a free monogenerated $LD$ monoid structure on a subset $\mathcal{F}_{e}$ of $R$.
Then in Section \ref{sec:order}, we extend the Dehornoy order from $B_{\infty}$ to the whole monoid $R$ by using the Laver order on $\mathcal{F}_{e}$ and the action of $R$ on the free group $F_{\infty}$.
We also prove this action to be faithful by showing that the Laver order on $\mathcal{F}_{e}$ is equivalent to an order defined by the action (Lemma \ref{lem:coincide-on-F_e}).
Finally, in Section \ref{sec:relations}, we prove that this action factors through the quotient map from $R$ to the monoid of shrinking braids and deduce Theorem \ref{thm:main}.
\subsection*{Acknowledgement}
This work was finished when the author was at Tsinghua University. The author thanks anonymous referees for pointing out related works by P. Dehornoy. 
\section{Shrinking braid monoid}\label{sec:shrink}

\begin{figure}
    \centering
    \includegraphics{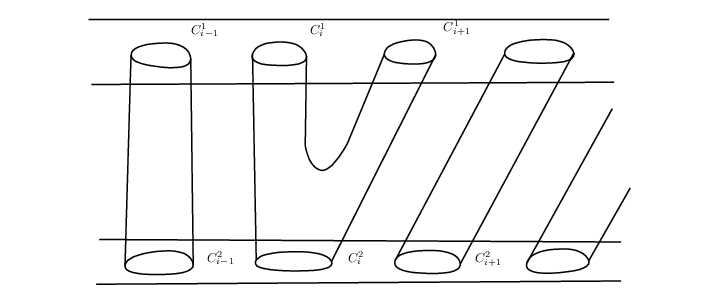}
    \caption{The cobordism $\tilde{x_{i}}$.}
    \label{fig:pant}
\end{figure}

\begin{figure}
    \centering
    \includegraphics{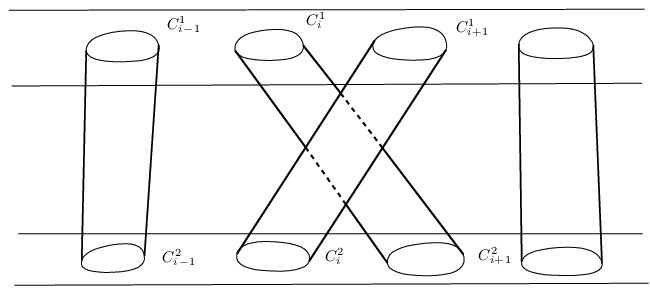}
    \caption{The cobordism $\tilde{s_{i}}$.}
    \label{fig:braid}
\end{figure}

We consider planes $\Gamma_{1}$ and $\Gamma_{2}$ in $\mathbb{R}^{3}$
where $\Gamma_{1} = \{ ( x,y,1 ) | x,y \in \mathbbm{R} \nobracket \}$ and
$\Gamma_{2} = \{ ( x,y,0 ) | x,y \in \mathbbm{R} \nobracket \}$. Denote by $\sum_{1}
= \{ C^{1}_{i} | i=1,2, \ldots \nobracket \}$ the set of input circles where
$C_{i}^{1}$ is the circle with radius $1/3$ centered at $( i,0,1 )$ in the
plane $\Gamma_{1}$. Denote by $\sum_{2} = \{ C^{2}_{i} | i=1,2, \ldots \nobracket
\}$ the set of output circles where $C_{i}^{2}$ is the circle with radius $1/3$
centered at $( i,0,0 )$ in the plane $\Gamma_{2}$.
Let $N= [ 0, \infty ) \times [ -1,1 ] \times [ 0,1 ]$.
\begin{definition}\label{def:cobordisms}
  We call $r$ a cobordism if we can write $r=\bigcup_{i=1}^{\infty}\rho_{i}$ where $\{\rho_{i}: i=1,2,\ldots\}$ is a countable collection of embedded smooth surfaces inside $N$ such that the following holds:
  \begin{enumerate}
      \item For each $i\in \{1,2,\ldots\}$, there are two finite subsets $I^{1}_{i}$ and $I^{2}_{i}$ of $\{1,2,\ldots\}$ such that $\partial \rho_{i}= (\bigcup_{k\in I^{1}_{i}}C^{1}_{k}) \bigcup (\bigcup_{j \in I^{2}_{i}}C^{2}_{j})$.
      \item $\rho_{i} \cap \rho_{i'}=\emptyset$ when $i\neq i'$.
      \item $\bigcup_{i=1}^{\infty} I^{1}_{i}= \bigcup_{i=1}^{\infty} I^{2}_{i}=\{1,2,\ldots\}$.
  \end{enumerate}
    We denote by $\mathcal{C}$ the set of all cobordisms.  
    Given two cobordisms $r_{1}$ and $r_{2}$, they are said to be equivalent if and only if, there exists a diffeomorphism $G:N \rightarrow N$ such that, $G(x)=x$ for any $x \in \partial N$ and $G(r_{1})=r_{2}$. We denote $r_{1}\sim r_{2}$ when $r_{1}$ and $r_{2}$ are equivalent. Let $\tilde{\mathcal{C}}=\mathcal{C}/\sim$ be the quotient of $\mathcal{C}$ under the equivalent relation $\sim$. Sometimes, we also call elements in $\tilde{\mathcal{C}}$ cobordisms if there is no confusion. 
    
    If $r_1,r_2 \in \tilde{\mathcal{C}}$, we let $r_1 r_2$ be the concatenation of $r_1$ and $r_2$. $\tilde{\mathcal{C}}$ is a monoid under this composition.
    
\end{definition}
We define two collections of elements in $\tilde{\mathcal{C}}$ which are important to us:

The cobordism $\tilde{x_{i}}$ (the \emph{$i$-th pant}) is defined as follow: for $k<i$
there is a straight tube from $C_{k}^{1}$ to $C_{k}^{2}$; for $k>i+1$ there is a
straight tube from $C_{k}^{1}$ to $C_{k-1}^{2}$; and there
is a pant-like cobordism from $C_{i}^{1}$ and $C_{i+1}^{1}$ to $C_{i}^{2}$. See Figure \ref{fig:pant} for an illustration.
The cobordism $\tilde{s_{i}}$ (the \emph{$i$-th braid}) is the braid-like tube cobordism twisting
between the $i$-th disk and $i+1$-th disk. See Figure \ref{fig:braid} for an illustration. And $\tilde{s_{i}}^{-1}$ is
defined as the reflection of $\tilde{s_{i}}$ along the plane $\{ (x,y,\frac{1}{2}) : x,y \in \R \}$. 

\begin{definition}
  Let $S \subseteq \tilde{\mathcal{C}}$ be the submonoid generated by $\{\tilde{x_{i}},\tilde{s_{i}},\tilde{s_{i}}^{-1}:i \in \mathbbm{N}\}$.   
\end{definition}
\begin{proposition}\label{prop:relationofS}
  The corbordisms $\tilde{x_{i}}$ and $\tilde{s_{i}}$ satisfy the following relations in $S$:
\begin{enumerate}
    \item $\tilde{s_{i}} \tilde{x_{i}} =\tilde{x_{i}}$\label{rel:1}
    \item $\widetilde{x_{i+1}} \tilde{x_{i}} = \tilde{x_{i}}^{2}$\label{rel:2}
    \item $\widetilde{x_{i+1}} \tilde{s_{i}} =\tilde{s_{i}} \widetilde{s_{i+1}} \tilde{x_{i}}$\label{rel:3}
    \item $\tilde{x_{i}} \tilde{s_{i}} =\widetilde{s_{i+1}} \tilde{s_{i}} \widetilde{x_{i+1}}$\label{rel:4}
    \item $\tilde{x_{i}} \tilde{x_{j}} =\widetilde{x_{j+1}} \tilde{x_{i}}$ , $\tilde{x_{i}} \tilde{s_{j}} =\widetilde{s_{j+1}} \tilde{x_{i}}$ when $j>i$ and $\tilde{x_{i}} \tilde{s_{j}} =\tilde{s_{j}} \tilde{x_{i}}$ for $j<i-1$\label{rel:5}
    \item $\tilde{s_{i}} \widetilde{s_{i+1}} \tilde{s_{i}} = \widetilde{s_{i+1}} \tilde{s_{i}} \widetilde{s_{i+1}}$\label{rel:6}
    \item $\tilde{s_{i}} \tilde{s_{j}} =\tilde{s_{j}} \tilde{s_{i}}$ $\tmop{for}  j>i+1$\label{rel:7}
\end{enumerate}
\end{proposition}
\begin{proof}
  Directly check by diagrams.
\end{proof}

\begin{definition}\label{def:define-R}
  The monoid $R$ is the monoid (algebraically) generated by $s_{i}$, $s_{i}^{-1}$ and
  $x_{i}$ $( i \in \mathbbm{N} )$ with the seven relations given above (by substituting $\tilde{x_{i}},\tilde{s_{i}}$ by $x_{i},s_{i}$ respectively) together with $s_{i}^{-1}s_{i}=s_{i} s_{i}^{-1}=e$ for each $i\in \mathbbm{N}$. Here, $e$ is the identity element in $R$.
  We also denote by $R^{+}$ the submonoid of $R$ generated by $s_{i}$ and $x_{i}$ $( i \in
  \mathbbm{N} )$.
\end{definition}
\begin{remark}\label{rem:relation-to-parenthesized}
In \cite{dehornoy2006group}, the parenthesized braid group $B_{\bullet}$ was constructed and its relations are specified in \cite[Lemma 1.15]{dehornoy2006group}. By using group generators $\{s_{i},x_{i}:i\in \Z_{+}\}$, $B_{\bullet}$ is represented by relations $(3)$ to $(7)$ in Proposition \ref{prop:relationofS}. The monoid $R$ is obtained from $B_{\bullet}$ by adding relation $(1)$ and $(2)$, and abandoning the invertibility of $x_{i}$'s. This is in agreement with the fact that the pant-like cobordisms $\tilde{x_{i}}$'s have no inverse in $S$.
\end{remark}
\begin{proposition}\label{prop:define-alpha}
  There is a surjective morphism $\alpha$ from $R$ to $S$ which maps $x_{i}$ to
  the $i$-th pant $\tilde{x_{i}}$ and maps $s_{i}$ to the $i$-th braid $\tilde{s_{i}}$ for each $i\geq 1$.
  \begin{proof}
    By Proposition \ref{prop:relationofS}.
  \end{proof}
\end{proposition}

In Section \ref{sec:relations}, we will prove the following theorem.
\begin{theorem}\label{thm:main}
  The monoid morphism $\alpha:R\rightarrow S$ defined in Proposition \ref{prop:define-alpha} is an isomorphism.
\end{theorem}

\section{Basic properties of $R$}\label{sec:basic}

In this section, ``morphism'' means a morphism between two monoids, i.e. a
function that preserves the multiplication and the identity. For example, there is a shift
morphism $\tmop{sh}:R \rightarrow R$ defined by $\tmop{sh} ( s_{i} ) =s_{i+1}
\nocomma$ and $\tmop{sh} ( x_{i} ) =x_{i+1}$.

We first investigate the submonoid of $R$ generated by $\{ x_{i} : i=1,2,
\ldots \nobracket \}$.

Let $X_{\infty}$ be the monoid defined by the representation 
\begin{equation}\label{eq:defX}
    X_{\infty}=\langle y_{i} | y_{i} y_{j} =y_{j+1} y_{i} , 1 \leq i \leq j <\infty \rangle.
\end{equation}
Then by relation $(2)$ and $(4)$ in Proposition \ref{prop:relationofS}, we have a morphism $\upsilon :X_{\infty}
\longrightarrow R$ such that $\upsilon(y_{i})=x_{i}$ for each $i \in \N$.

Let $F$ be the set of all the surjections from $\N$ to itself. Endow $F$ a monoid structure by defining the composition as follows
\begin{equation}
    f g := g\circ f.
\end{equation}
Then the identity element of $F$ is the identity function. For each $i \in \N$, define functions $g_{i},h_{i}\in F$ as follows:
\begin{enumerate}
    \item $g_{i} ( k ) =k$ for $k \leqslant i$ and $g_{i} ( k ) =k-1$ for $k>i$. 
    \item $h_{i}(i)=i+1$, $h_{i}(i+1)=i$ and $h_{i}(k)=k$ if $k \neq \text{$i$ or $i+1$}$.
\end{enumerate}
\begin{lemma}\label{lem:relationofF}
  The seven relations in Proposition \ref{prop:relationofS} hold in $F$ after substituting $\tilde{x_{i}}$ by $g_{i}$ and $\tilde{s_{i}}$ by $h_{i}$ for each $i\in \Z_{+}$.
\end{lemma}
\begin{proof}
  Easily check.
\end{proof}
By Lemma \ref{lem:relationofF}, there is a monoid morphism $\psi :R \longrightarrow F$ such that $\psi ( x_{i} ) =g_{i}$ and $\psi ( s_{i} )
=h_{i}$ for each $i \in \N$.
Let $\phi = \psi \circ \upsilon$ be the morphism from $X_{\infty}$ to $F$.
By definition, $\phi ( y_{i} ) =g_{i}$.

Let $y \in X_{\infty}$ and write $y=y_{i_{1}}y_{i_{2}}\cdots y_{i_{m}}$. For any $k=1,2,\ldots,m-1$, if $i_{k}>i_{k+1}$, by 
\eqref{eq:defX}, we can substitute $y_{i_{k}}y_{i_{k+1}}$ by $y_{i_{k+1}}y_{i_{k}-1}$. By iteration, we can write
\begin{equation}\label{eq:word-form-increasing-indice}
    y=y_{i'_{1}} y_{i'_{2}} \ldots y_{i'_{m}}
\end{equation}
where $i'_{1}
\leqslant i'_{2} \leqslant   \ldots \leqslant i'_{m}$.
Let $f= \phi ( y )=g_{i'_{1}}g_{i'_{2}}\cdots g_{i'_{m}}=g_{i'_{m}}\circ g_{i'_{m-1}}\circ \cdots \circ g_{i'_{1}}$.
Under this word form, for each $j \in \N$, we let $n'_{j}=\#\{k:i'_{k}=j\}$ which is a nonnegative integer. By denoting $y_{i}^{0}=e$, we can rewrite
\begin{equation}\label{eq:wordform}
    y=y_{1}^{n'_{1}} y_{2}^{n'_{2}}\cdots y_{M}^{n'_{M}}
\end{equation}
for some $M \in \N$ and thus
\begin{equation}\label{eq:def-f}
    f=g_{M}^{n'_{M}}\circ \cdots \circ g_{2}^{n'_{2}}\circ g_{1}^{n'_{1}}.
\end{equation}

\begin{lemma}\label{lem:f's image}
Let $f:\N \rightarrow \N$ be the function in \eqref{eq:def-f} with nonnegative integers $n'_{1},\ldots,n'_{M}$. Further define $n'_{i}=0$ for $i\geq M+1$. Then we have
$$f ( m ) = \tmop{min} \{ l\in \N : n'_{1} +n'_{2} + \ldots
+n'_{l}+l \geqslant m  \}$$ for any $m\in \N$.
\end{lemma}
\begin{proof}
  By definition, for any $m,i\in \N$ and nonnegative integer $p$, we have
  \begin{equation}\label{eq:g-acting}
      g_{i}^{p}(m)=
      \begin{cases}
      \max\{m-p,i\}& \text{\quad if $m\geq i$}\\ 
      m& \text{\quad if $m < i$}. 
      \end{cases}
  \end{equation}
   Suppose $l \in \N$ satisfy 
  \begin{equation}\label{eq:lmin1}
        n'_{1} +n'_{2} + \ldots +n'_{l}+l \geqslant m
  \end{equation}
and
\begin{equation}\label{eq:lmin2}
    n'_{1} +n'_{2} + \ldots +n'_{l-1}+l-1 < m.
\end{equation}
Then using \eqref{eq:g-acting} iteratively, we have 
\begin{equation}
    g_{l-1}^{n'_{l-1}} \circ g_{l-2}^{n'_{l-2}} \cdots \circ g_{1}^{n'_{1}}(m)=m-n'_{1}-n'_{2}-\ldots-n'_{l-1} \geq l
\end{equation}
where the second inequality is from \eqref{eq:lmin2}. Denote $m'=m-n'_{1}-n'_{2}-\ldots-n'_{l-1}$. By \eqref{eq:g-acting} and \eqref{eq:lmin1}, since $m' \geq l$, we have
$g_{l}^{n'_{l}}(m')=\max\{m'-n'_{l},l\}=l$. Hence 
$$f(m)=g_{M}^{n'_{M}}\circ \cdots \circ g_{2}^{n'_{2}}\circ g_{1}^{n'_{1}}(m)=g_{M}^{n'_{M}}\circ \cdots \circ g_{l}^{n'_{l}}(m')=g_{M}^{n'_{M}}\circ \cdots \circ g_{l+1}^{n'_{l+1}}(l)=l$$ 
where the last equality is due to \eqref{eq:g-acting}.
Our lemma follows.
\end{proof}
\begin{lemma}\label{lem:n' value}
Let $f:\N \rightarrow \N$ be the function in \eqref{eq:def-f} with nonnegative integers $n'_{1},\ldots,n'_{M}$. Further define $n'_{i}=0$ for $i\geq M+1$. Then we have
$$\#f^{-1}(k)=n'_{k}+1$$ for $k\in \N$.
\end{lemma}
\begin{proof}
By Lemma \ref{lem:f's image}, we have
  $$f^{-1}(k)=\{n'_1+\ldots +n'_{k-1}+k,n'_1+\ldots +n'_{k-1}+k+1,\ldots,n'_1+\ldots +n'_{k}+k \}.$$ Our lemma follows.
\end{proof}
\begin{proposition}
  The word form \eqref{eq:wordform} is unique for any $y \in X_{\infty}$ and $\phi$ is
  injective.
\end{proposition}
  \begin{proof}
    Let $f= \phi ( y )$ and then $n'_{i} =\#f^{-1} ( i )-1 $ by Lemma \ref{lem:n' value}. Thus the word form \eqref{eq:wordform} of any element $y$ in $X_{\infty}$ is
    determined by $\phi(y)$ and thus is unique. For any $y_{1},y_{2}\in X_{\infty}$, if $\phi(y_{1})=\phi(y_{2})$, then the word forms of $y_{1}$ and $y_{2}$ are the same and thus $y_{1}=y_{2}$. This implies $\phi$ is injective.
  \end{proof}

Because $\phi$ is injective and $\phi = \psi \circ \upsilon$, we deduce that $\upsilon$ is
injective and $X_{\infty}$ is embedded in $R$ as a sub-monoid. The following theorem says that the same is true for $B_{\infty}$.

\begin{theorem}\label{thm:B_inf-embed}
  $B_{\infty}$ is embedded in $R$. i.e. The monoid generated by $\{ s_{i},s_{i}^{-1} : i=1,2 \ldots\}$ in $R$ is isometric to $B_{\infty}$ (as a monoid). Moreover, suppose $r\in R$ and $\Bar{r}_{1},\Bar{r}_{2}$ are two word forms of $r$ in $R$. Then $\Bar{r}_{1}$ and $\Bar{r}_{2}$ have the same total number of $x_{i}$'s.  
  
\end{theorem}
\begin{proof}
  Let $B_{\infty}'\subset R$ be the monoid generated by $s_{i},s_{i}^{-1}$ for $i\in \N$.
    Suppose $x,y\in B_{\infty}'$ and they are equivalent in
    $R$. Then there is a procedure from $x$ to $y$ by using the rules $( 1 ) \sim ( 7 )$ in Proposition \ref{prop:relationofS} plus the rule $s_{i}
    s_{i}^{-1} =s_{i}^{-1} s_{i} =e$.

    Observe that the rules $( 1 ) \sim ( 5 )$ never show up in the procedure since the words on
    both sides of the rule contain elements in $X_{\infty}$. This means that, from
    $x$ to $y$, the procedure only uses rules $( 6 ) , ( 7 )$ and the rule $s_{i} s_{i}^{-1}
    =s_{i}^{-1} s_{i} =e.$ Thus, by definition of the braid group $B_{\infty}$, we have $B_{\infty}' \cong B_{\infty}$.

    Finally note that under any of the seven rules in Proposition \ref{prop:relationofS}, the total number of $x_{i}$'s does not change. Our theorem follows.
  \end{proof}

\section{Enveloping LD monoid}\label{sec:env}

A {\tmem{left-distributive monoid(LD monoid)}} is a set with two compositions
$\circ$ and $\cdot$ such that they satisfy the laws:
\begin{align}
    a \cdot ( b \circ c ) =& ( a \cdot b ) \circ ( a \cdot c )\label{eq:law(2)} \\
    ( a \cdot b ) \circ a =& a \circ b\label{eq:law(3)} \\
    a  \cdot ( b \cdot c ) =& ( a \circ b )\cdot c\label{eq:law(4)}
\end{align}
and the associative law of $\circ$.

From the three laws above can we imply the \emph{left distributive law} of $\cdot$ :
\begin{equation}\label{eq:law(1)}
  a \cdot ( b \cdot c ) = ( a \cdot b ) \cdot ( a \cdot c )
\end{equation}
(proof: $a \cdot ( b \cdot c ) = ( a \circ b ) \cdot c= ( ( a \cdot b ) \circ
a ) \cdot c= ( a \cdot b ) \cdot ( a \cdot c )$)

A set with a single $\cdot$ composition satisfying \eqref{eq:law(1)} is called a
{\tmem{left-distributive system}} (or \emph{LD system}).

Let $( \mathcal{A}, \cdot )$ be an LD system, we will
construct an LD monoid $\mathcal{F}$ generated by $\mathcal{A}$ which is
universal in the following sense. For any LD monoid $\mathcal{G}$ and any LD system
homomorphism $f:\mathcal{A} \rightarrow \mathcal{G}$, there is a unique LD monoid
homomorphism $g:\mathcal{F} \rightarrow \mathcal{G}$ such that
$f=g \circ i$ where $i$ is the canonical embedding of $\mathcal{A}$ into
$\mathcal{F}$.

Let $\mathfrak{F}$ be the free semigroup generated by the elements
of $\mathcal{A}$ with composition denoted by $\circ$. Then we can write $\mathfrak{F}= \bigcup_{i=1}^{\infty}
\mathcal{A}^{\circ i}$ where $\mathcal{A}^{\circ n} = \{ a_{1} \circ a_{2}
\circ   \ldots \circ a_{n} | a_{i} \in \mathcal{A} \nobracket \}$. Now define
the composition $\cdot$ on $\mathfrak{F}$ by
\begin{align}\label{eq:def-composition}\tag{$\diamondsuit$}
\begin{split}
          \ &( a_{1} \circ  a_{2} \circ   \ldots \circ a_{n} ) \cdot ( b_{1} \circ b_{2}
   \circ   \ldots b_{m} ) =\\
    &( a_{1} \cdot ( a_{2} \cdot \ldots ( a_{n} \cdot
   b_{1} ) ) ) \circ ( a_{1} \cdot ( a_{2} \cdot \ldots \cdot ( a_{n} \cdot
   b_{2} )))   \circ   \ldots \circ ( a_{1} \cdot ( a_{2} \cdot \ldots \cdot (
   a_{n} \cdot b_{m} ) ) ) \  
\end{split}
\end{align}



Consider the positive braid monoid $B_{n}^{+}$ generated by $s_{i} (i<n)$ with relations $s_{i} s_{i+1} s_{i}=s_{i+1}s_{i}s_{i+1}  $ and $s_{i}s_{j}=s_{j}s_{i}$ for $|i-j|>1$. Define an action of $B_{n}^{+}$ on
$\mathcal{A}^{\circ n}$ by
\begin{equation}\label{eq:action-B_n+}
     s_{i} ( a_{1} \circ a_{2} \circ   \ldots \circ a_{i} \circ a_{i+1}
   \circ   \ldots \circ a_{n} ) =a_{1} \circ a_{2} \circ   \ldots \circ a_{i}
   \cdot a_{i+1} \circ a_{i} \circ   \ldots \circ a_{n}. 
\end{equation}
By the left distributive law of $\mathcal{A}$, one can check that
\begin{align}
    s_{i} (s_{i+1} (s_{i} (\mathfrak{a})))=& s_{i+1}(s_{i}(s_{i+1}(\mathfrak{a}))) &\text{ for $1\leq i\leq n-2$}\\
    s_{i} (s_{j} (\mathfrak{a}))=&
    s_{j} (s_{i} (\mathfrak{a})) &\text{ for $1\leq i,j\leq n-1$ with $|i-j|>1$ }
\end{align}
for any $\mathfrak{a}\in \mathcal{A}^{\circ n}$. Thus the group action given by \eqref{eq:action-B_n+} is well-defined. Now we
define an equivalent relation on $\mathcal{A}^{\circ n}$ for each $n\in \N$. Given $\mathfrak{b}_{1},\mathfrak{b}_{2}\in \mathcal{A}^{\circ n}$, define $\mathfrak{b}_{1}\equiv \mathfrak{b}_{2}$ if and only if there are $g_{1},g_{2}\in B_{n}^{+}$ such that $g_{1}(\mathfrak{b}_{1})=g_{2}(\mathfrak{b}_{2})$.
Let $\widetilde{\mathcal{A}^{\circ n}}$ be the set of equivalent classes of
$\mathcal{A}^{\circ n}$ under the equivalent relation $\equiv$, and let $\mathcal{F}=
\bigcup_{i=1}^{\infty} \widetilde{\mathcal{A}^{\circ i}}$. Note that $\widetilde{\mathcal{A}^{\circ 1}} =\mathcal{A}$,
thus $\mathcal{A}$ is embedded into $\mathcal{F}$. Endow the compositions $\cdot$ and
$\circ$ on $\mathcal{F}$ as follows:
\begin{align}
      \overline{\mathfrak{b}_{1}} \cdot \overline{\mathfrak{b}_{2}}=& \overline{\mathfrak{b}_{1} \cdot \mathfrak{b}_{2}} \label{eq:quotient-1} \\
      \overline{\mathfrak{b}_{1}} \circ \overline{\mathfrak{b}_{2}}=& \overline{\mathfrak{b}_{1} \circ \mathfrak{b}_{2}}.\label{eq:quotient-2}
\end{align}
Here, $\mathfrak{b}_{1},\mathfrak{b}_{2}\in \mathfrak{F}$ and $\overline{\mathfrak{b}_{i}}\in \mathcal{F}$ is the equivalent class of $\mathfrak{b}_{i}$ for $i=1,2$. 
\begin{proposition}\label{prop:equal-in-F}
  Let $m$ be a positive integer and $a_{i} ,b_{i}$ be elements in $\mathcal{A}$ for $1\leq i\leq m$. Then $\overline{a_{1} \circ a_{2}
  \circ   \ldots \circ a_{m}} =\overline{b_{1} \circ b_{2} \circ   \ldots \circ b_{m}}$ if and only if there are $g_{1} ,g_{2} \in B_{m}^{+}$
  such that $g_{1} ( a_{1} \circ a_{2}
  \circ   \ldots \circ a_{m} ) =g_{2} ( b_{1} \circ b_{2} \circ   \ldots \circ b_{m} )$.

  \begin{proof}
    By definition of equivalence relation $\equiv$.
  \end{proof}
\end{proposition}
\begin{theorem}
  Suppose $( \mathcal{A}, \cdot )$ is an LD system and $\mathcal{F}$ is
  constructed as in the above argument. Then the compositions $\cdot$ and
  $\circ$ are well-defined on $\mathcal{F}$ and $( \mathcal{F}, \cdot , \circ
  )$ is an LD monoid.

\end{theorem}
\begin{proof}
    The well-definedness of $\circ$ by \eqref{eq:quotient-2} is evident and we now check that $\cdot$
    is well defined.
    \begin{enumerateroman}
      \item Given $a_{i} ,b \in \mathcal{A}$

      \ \ \ \ \ \ \ \ $s_{i} ( a_{1} \circ a_{2} \circ   \ldots \circ
      a_{i} \circ a_{i+1} \circ   \ldots \circ a_{n} ) \cdot b$

      \ \ \ \ \ \ \ $= ( a_{1} \circ a_{2} \circ   \ldots \circ a_{i} \cdot
      a_{i+1} \circ a_{i} \circ   \ldots \circ a_{n} ) \cdot b$

      \ \ \ \ \ \ \ $= a_{1} \cdot ( a_{2} \cdot ( \ldots \cdot ( ( a_{i}
      \cdot a_{i+1} ) \cdot ( a_{i} \cdot \ldots \cdot ( a_{n} \cdot b )   )
      )   )   )$ \ \ \ \ by $( \diamondsuit )$

      \ \ \ \ \ \ \ $=  a_{1} \cdot ( a_{2} \cdot ( \ldots \cdot ( a_{i}
      \cdot ( a_{i+1} \cdot \ldots \cdot ( a_{n} \cdot b )   )   )   )
      \nobracket$ \ \ \ \ \ \ \ \ \ \ \ by \eqref{eq:law(1)} on $\mathcal{A}$

      \ \ \ \ \ \ \ $=  ( a_{1} \circ a_{2} \circ   \ldots \circ a_{i} \circ
      a_{i+1} \circ   \ldots \circ a_{n} ) \cdot b$ \ \ \ \ \ \ \ \ \ \ \ \ \
      \ \ \ \ by $( \diamondsuit )$

      \item Given $b_{i} \in \mathcal{A}$ and $a \in \mathfrak{F}$

      $a \cdot s_{i} ( b_{1} \circ b_{2} \circ   \ldots \circ b_{i} \circ
      b_{i+1} \circ   \ldots \circ b_{n} )$

      $=a \cdot ( b_{1} \circ b_{2} \circ   \ldots \circ b_{i} \cdot b_{i+1}
      \circ b_{i} \circ   \ldots \circ b_{n} )$

      $=  ( a \cdot b_{1} \circ a \cdot b_{2} \circ   \ldots \circ a \cdot (
      b_{i} \cdot b_{i+1} ) \circ a \cdot b_{i} \circ   \ldots \circ a \cdot
      b_{n} )$ \ \ \ \ \ \ \ \ \ \ \ \ \ \ \ \ \ by $( \diamondsuit )$

      $= ( a \cdot b_{1} \circ a \cdot b_{2} \circ   \ldots \circ ( \nobracket
      ( \nobracket a \cdot b_{i} )   \cdot ( \nobracket a \cdot b_{i+1} )
      \nobracket )   \circ a \cdot b_{i} \circ   \ldots \circ a \cdot b_{n} )$
      \ \ \ \ \ \ \ \ \ \ \ by \eqref{eq:law(1)} which can be deduced on $\mathfrak{F}$

      $= s_{i} ( a \cdot b_{1} \circ a \cdot b_{2} \circ   \ldots \circ (
      \nobracket a \cdot b_{i} )   \circ ( \nobracket a \cdot b_{i+1} ) \circ
      \ldots \circ a \cdot b_{n} )$
      
      $= s_{i} ( a \cdot ( b_{1} \circ b_{2} \circ   \ldots \circ b_{i} \circ
      b_{i+1} \circ   \ldots \circ b_{n} ))$
    \end{enumerateroman}
    From $i$ and $ii$, we conclude that the composition $\cdot$ is well
    defined by \eqref{eq:quotient-1} on $\mathcal{F}$.

    Now we prove that $( \mathcal{F}, \cdot , \circ )$ is an LD monoid.
    Firstly, \eqref{eq:law(2)} and \eqref{eq:law(4)} hold because of definition of composition (equation \eqref{eq:def-composition}) on
    $\mathfrak{F}$ and the definition of $\mathcal{F}$. We now check \eqref{eq:law(3)} holds.

    Write $\mathfrak{a}=a_{1} \circ a_{2} \circ    \ldots \circ a_{n}$ and $\mathfrak{b}=b_{1} \circ
    b_{2} \circ   \ldots \circ b_{m}$. Then we have \ \ \

    $( \mathfrak{a} \cdot \mathfrak{b} ) \circ \mathfrak{a}= ( a_{1} \cdot ( a_{2} \cdot \ldots ( a_{n} \cdot
    b_{1} ) ) ) \circ ( a_{1} \cdot ( a_{2} \cdot \ldots \cdot ( a_{n} \cdot
    b_{2} )))   \circ   \ldots \circ ( a_{1} \cdot ( a_{2} \cdot \ldots \cdot (
    a_{n} \cdot b_{m} ) ) ) \circ a_{1} \circ a_{2} \circ   \ldots \circ
    a_{n}$

    $=  ( a_{1} \cdot ( a_{2} \cdot \ldots ( a_{n} \cdot b_{1} ) ) ) \circ (
    a_{1} \cdot ( a_{2} \cdot \ldots \cdot ( a_{n} \cdot b_{2} )))   \circ
    \ldots \circ (a_{1} \circ ( a_{2} \cdot \ldots \cdot ( a_{n} \cdot b_{m} ))
    ) \circ a_{2} \circ   \ldots \circ a_{n}$ \ \ (This is by definition of $\mathcal{F}$, i.e. the equivalence relation $\equiv$)

    $=  ( a_{1} \cdot ( a_{2} \cdot \ldots ( a_{n} \cdot b_{1} ) ) ) \circ (
    a_{1} \cdot ( a_{2} \cdot \ldots \cdot ( a_{n} \cdot b_{2} ) ) )  \circ
    \ldots \circ a_{1} \circ (a_{2} \circ ( a_{3} \cdot \ldots \cdot ( a_{n}
    \cdot b_{m} ) )) \circ a_{3} \circ   \ldots \circ a_{n}$

    $= \ldots$

    $= ( a_{1} \cdot ( a_{2} \cdot \ldots ( a_{n} \cdot b_{1} ) ) ) \circ (
    a_{1} \cdot ( a_{2} \cdot \ldots \cdot ( a_{n} \cdot b_{2} ) ) )  \circ
    \ldots \circ ( a_{1} \cdot ( \nobracket a_{2} \ldots \cdot ( a_{n} \cdot
    b_{m-1} ) ) \nobracket ) \circ a_{1} \circ a_{2} \circ   \ldots \circ
    a_{n}\circ b_{m}$

    $= \ldots$

    $=a_{1} \circ a_{2} \circ   \ldots \circ a_{n} \circ b_{1} \circ b_{2}
    \circ   \ldots b_{m}$

    $=\mathfrak{a} \circ \mathfrak{b}$

    Thus \eqref{eq:law(3)} holds and $( \mathcal{F}, \cdot , \circ )$ is an LD monoid. 
  \end{proof}
We call the LD monoid \tmverbatim{$\mathcal{F}$} constructed as in the above
argument the \emph{enveloping LD monoid} of the LD system $\mathcal{A}$.

\begin{theorem}\label{thm:universal}
  Suppose $\mathcal{A}$ is an LD system and $\mathcal{F}$ is the enveloping LD monoid generated by $\mathcal{A}$. Then for any LD monoid $\mathcal{G}$ and any
  LD system homomorphism $f:\mathcal{A} \rightarrow \mathcal{G}$, there is a unique
  homomorphism of LD monoid $g:\mathcal{F} \rightarrow \mathcal{G}$ such that
  $f=g \circ i$ where $i$ is the canonical embedding of $\mathcal{A}$ into
  $\mathcal{F}$. 
\end{theorem}
\begin{proof}
    The uniqueness of $g$ is evident and we will prove its existence.
    Define the function $g$ by $g (\overline{ a_{1} \circ a_{2} \circ    \ldots \circ
    a_{n} }) =f ( a_{1} ) \circ f ( a_{2} ) \circ   \ldots \circ f ( a_{n} )$
    for any $a_{i} \in \mathcal{A}$. This is well defined because
    \begin{align}
        \begin{split}
           &g (\overline{ a_{1} \circ a_{2} \circ   \ldots \circ
    a_{i} \cdot a_{i+1} \circ a_{i} \circ   \ldots \circ a_{n} })\\
   =&f ( a_{1} ) \circ f ( a_{2} ) \circ   \ldots \circ  f
    ( a_{i} \cdot a_{i+1} ) \circ  f ( a_{i} ) \circ \ldots \circ f ( a_{n} )\\
     =&f ( a_{1} ) \circ f ( a_{2} ) \circ   \ldots \circ   (
    f ( a_{i} )   \cdot f ( a_{i+1} ) ) \circ  f ( a_{i} ) \circ \ldots \circ
    f ( a_{n} ) \\
   =&f ( a_{1} ) \circ f ( a_{2} ) \circ   \ldots \circ  f (
    a_{i} )   \circ  f ( a_{i+1} ) \circ \ldots \circ f ( a_{n} )\\
   =  &g (\overline{ a_{1} \circ a_{2} \circ   \ldots \circ a_{i}
    \circ a_{i+1} \circ   \ldots \circ a_{n} }) 
        \end{split}
    \end{align}
    where the third equation is because $\mathcal{G}$ is an LD monoid. Since $f=g \circ i$, the conclusion follows.
\end{proof}
\begin{corollary}
  Suppose $\mathcal{A}$ is the monogenerated free LD system and $\mathcal{F}$ is the enveloping LD monoid generated by $\mathcal{A}$. Then
  $\mathcal{F}$ is the monogenerated free LD monoid.
\end{corollary}
\begin{proof}
    Suppose $\mathcal{A}$ is a monogenerated free LD system. Then there is a monogenerated LD system $\mathcal{A}'$
    which is embedded in a monogenerated free LD monoid $\mathcal{F}'$ with
    the same generator (\cite{laver1992left}). Consider the LD system isomorphism $\phi$ from $\mathcal{A}$ to
    $\mathcal{A}'$ which maps the generator of $\mathcal{A}$ to the generator of $\mathcal{A}'$. By Theorem \ref{thm:universal}, there is an LD monoid homomorphism $\Phi$ from
    $\mathcal{F}$ to $\mathcal{F}'$ which extends $\phi$.

    Since $\mathcal{F}'$ is free, there is an LD monoid homomorphism $\Phi':\mathcal{F}'\rightarrow \mathcal{F}$ which maps the generator of $\mathcal{F}'$ to the generator of $\mathcal{F}$. Then $\Phi' \circ \phi$ is the identity map on $\mathcal{A}$. For any $b\in \mathcal{F}$, write $b=a_{1}\circ \cdots \circ a_{n}$ with $a_{i}\in \mathcal{A}$. Then
    \begin{equation}
        \Phi' \circ \Phi(b)=(\Phi' \circ \phi(a_{1}))\circ \cdots \circ (\Phi' \circ \phi(a_{n}))=a_{1}\circ\cdots \circ a_{n}=b.
    \end{equation}
    Thus $\Phi' \circ \Phi$ is the identity map and $\Phi$ is a bijection.
\end{proof}

\section{The LD monoid in $R$}\label{sec:monoidR}

In this chapter, we construct a subset of $R$ and endow it with the free monogenerated LD monoid
structure. An important result about the free monogenerated $\tmop{LD}$ monoid is the
existence of the Laver order:

\begin{theorem}[\cite{dehornoy2006group}\cite{laver1992left}]\label{thm:linear-Laver}
  Let $\mathcal{F}$ be the free monogenerated LD monoid. For any $p,q
  \in$\tmverbatim{$\mathcal{F}$}, define $p<q$ if and only if there are $p_{1}
  , \ldots p_{k} \in \mathcal{F}$ and $a_{1},\ldots,a_{k-1} \in \mathcal{F}$ such that following holds:
  $p_{1}=p$, $p_{k}=q$ and $p_{i+1}= p_{i} \cdot a_{i}$ for $1\leq i\leq k-2$, and $p_{k}=p_{k-1} * a_{k-1}$ with $*=\cdot \text{ or } \circ$.
  
  Then $<$ is a linear order on $\mathcal{F}$.
\end{theorem}

On the braid group
\begin{align*}
    &B_{\infty}\\
    =&\left< s_{1},s_{2},\ldots \big| s_{i} s_{i+1} s_{i}= s_{i+1} s_{i} s_{i+1}, s_{i} s_{j} =s_{j} s_{i} \text{ for all $i,j\geq 1$ with $|i-j|>1$}\right>,
\end{align*}
a left-distributive composition $\cdot$ has
been constructed (\cite{dehornoy1994braid}):
\begin{equation}\label{eq:composition-cdot}
     a \cdot b=a \tmop{sh} ( b ) s_{1} \tmop{sh} ( a^{-1} ).
\end{equation}
Denote the subset of $B_{\infty}$ generated by the identity under composition
$\cdot$ by $\mathcal{A}_{e}$. It was proved that $\mathcal{A}_{e}$ is a free
monogenerated LD system (\cite{dehornoy1994braid}).

The Laver order restricted on the monogenerated LD system $\mathcal{A}$ has a
braid group ordering counterpart:

\begin{theorem}[\cite{dehornoy1994braid}]\label{thm:linear-Dehornoy}
  There is a linear order $<_{D}$ on $B_{\infty}$ which is invariant under the left multiplication of $B_{\infty}$. Under this order, $a<_{D} b$ if and only if for some $k\in \mathbb{N}$, $a^{-1} b$ can be
  represented using words $s_{k} ,s_{k+1}^{\pm 1} ,s_{k+2}^{\pm 1} , \ldots$ and
  $s_{k}$ shows up in the word form.
\end{theorem}

If we consider the automorphism of $B_{n}$ which sends $s_{i}$ to $s_{n-i}$,
we imply the following:

\begin{theorem}[\cite{dehornoy1994braid}]\label{thm:word-form-in-B_inf}
For any element $b \in B_{\infty}$
which is not the identity, there exists $k\in \Z_{+}$ such that $b$ can be written in the word form of $s_{1}^{\pm}
,s_{2}^{\pm} , \ldots ,s_{k}^{\pm}$ where only one of $s_{k}$ and
$s_{k^{}}^{-1}$ appears in the form.
\end{theorem}

For each $n\in \Z_{+}$, consider the set of sequences of braids $$\mathfrak{B}^{(n)}=\{\vec{a}=( a_{1} ,a_{2} , \ldots ,a_{n} ):a_{i} \in B_{\infty}\}.$$ Define function $\tmop{sh}:\mathfrak{B}^{(n)} \rightarrow B_{\infty}$ by 
\begin{equation}\label{eq:def-sh}
    \tmop{sh} ( \vec{a} ) =a_{1} \tmop{sh} ( a_{2} ) \ldots
  \tmop{sh}^{n-1} ( a_{n} ).
\end{equation}
We define an action of $B_{n}^{+}$ on $\mathfrak{B}^{(n)}$ as follows (see \cite[Lemma 2.8]{dehornoy2018braid}),
\begin{equation}\label{eq:action-of-braid}
    s_{i} ( a_{1} ,a_{2} , \ldots ,a_{i} ,a_{i+1} , \ldots ,a_{n} ) = ( a_{1} ,
\ldots ,a_{i} \cdot a_{i+1} ,a_{i} , \ldots ,a_{n} )
\end{equation}
for $i<n$. The well-definedness follows from the left-distributive property of $\cdot$ composition. 

\begin{lemma}[\cite{dehornoy1994braid}]\label{lem:uniqueness-of-sh}
  For a sequence of braids $\vec{a} = ( a_{1} ,a_{2} , \ldots ,a_{n} )$,
  we have $\tmop{sh} ( g ( \vec{a} ) ) =
  \tmop{sh} ( \vec{a} ) g$ for any $g \in B_{n}^{+}$. Moreover, if $\vec{a},\vec{b} \in \mathfrak{B}^{(n)}$ with $a_{i},b_{i}\in \mathcal{A}_{e}$ for each $1\leq i\leq n$ and $\tmop{sh} ( \vec{a} ) = \tmop{sh} ( \vec{b}
  )$, then $\vec{a} = \vec{b}$.

\end{lemma}
\begin{proof}
    For the first property, see \cite[Lemma 2.10]{dehornoy2018braid}. For the second one, by \cite[Lemma 3.5]{dehornoy2018braid}, the Laver order on $\mathcal{A}_{e}$
    equals the $s_{1}$-positive order of $B_{\infty}$ restricting on
    $\mathcal{A}_{e}$ defined by $c_{1} <_{s_{1}} c_{2}$ if and only if there is a word form of $c_{1}^{-1} c_{2}$ containing
    $s_{1}$ but not containing $s_{1}^{-1}$. If $\vec{a} \not= \vec{b}$, let $1\leq i'<n$ be the minimal integer such that $a_{i'}\not=b_{i'}$. By Theorem \ref{thm:linear-Laver}, we can assume without loss of generality that $a_{i'}<_{s_{1}} b_{i'}$, i.e. a word form of $a_{i'}^{-1}b_{i'}$ contains $s_{1}$ but does not contain $s_{1}^{-1}$. However, this implies a word form of $(\tmop{sh} ( \vec{a} ))^{-1} \tmop{sh} ( \vec{b} )$ contains $s_{i'}$ but does not contain $s_{1}^{\pm},\ldots,s_{i'-1}^{\pm},s_{i'}^{-1}$. By Theorem \ref{thm:linear-Dehornoy},  $\tmop{sh} ( \vec{a} ) <_{D} \tmop{sh} ( \vec{b} )$.
  \end{proof}
We now construct an associative composition $\circ$ on monoid $R$ (defined in Definition \ref{def:define-R}). Recall that by Theorem \ref{thm:B_inf-embed}, $B_{\infty}$ is naturally embedded in $R$. The composition $\circ$ will satisfy the following properties: for any $a,b,c\in B_{\infty}\subset R$,

$( i )   ( a \cdot b ) \circ a=a \circ b $

$ ( ii )  a \cdot ( b \circ c ) = ( a \cdot b ) \circ ( a \cdot c )$.

Once this is done, we use equation
\eqref{eq:def-composition} to extend $\cdot$ to elements in $R$ with the form $a_{1}
\circ a_{2} \circ \ldots \circ a_{n}$ where $a_{i} \in B_{\infty}$.

Now, the composition $\circ$ is explicitly given by
\begin{equation}\label{eq:def-compo-circ}
    a \circ b=a \tmop{sh} ( b ) x_{1}.
\end{equation}
We now check that $\circ$ is associative and $( i ) , ( ii )$ hold for
$a,b,c \in B_{\infty}$.

\subparagraph{Associativity:}Given $a,b,c \in R$,

$( a \circ b ) \circ c=a \tmop{sh} ( b ) x_{1} \tmop{sh} ( c ) x_{1} =a
\tmop{sh} ( b ) \tmop{sh}^{2} ( c ) x^{2}_{1}$ \ \ \ \ \ \ \ \ (by relation \eqref{rel:5})

$=a \tmop{sh} ( b ) \tmop{sh}^{2} ( c ) x_{2} x_{1}$ \ \ \ \ \ \ \ \ \ \ \ \ \
\ \ \ \ \ \ \ \ \ \ \ \ \ \ \ \ \ \ \ \ \ \ \ \ \ \ \ \  (by relation \eqref{rel:2})

=$a \circ ( b \circ c )$

\subparagraph{(i):}Given $a,b \in B_{\infty}$,

$( a \cdot b ) \circ a=a \tmop{sh} ( b ) s_{1} \tmop{sh} ( a^{-1} ) \tmop{sh}
( a ) x_{1}$

$=a \tmop{sh} ( b ) s_{1} x_{1}$

$=a \tmop{sh} ( b ) x_{1}$ \ \ \ \ \ \ \ \ \ \ \ \ \ \ \ \ \ \ \ \ \ \ \ \ \ \
\ \ \ \ \ \ \ \ \ \ \ \ \ \ \ (by relation \eqref{rel:1})

$= a \circ b$

\subparagraph{(ii):}Given $a,b,c \in B_{\infty}$

$a \cdot ( b \circ c ) =a \tmop{sh} ( b \tmop{sh} ( c ) x_{1} ) s_{1}
\tmop{sh} ( a^{-1} )$

$=a \tmop{sh} ( b ) \tmop{sh}^{2} ( c ) x_{2} s_{1} \tmop{sh} ( a^{-1} )$

$=a \tmop{sh} ( b ) \tmop{sh}^{2} ( c ) s_{1} s_{2} x_{1} \tmop{sh} ( a^{-1}
)$ \ \ \ \ \ \ \ \ \ \ \ \ \ \ \ (by relation \eqref{rel:3})

$=a \tmop{sh} ( b ) \tmop{sh}^{2} ( c ) s_{1} s_{2} \tmop{sh}^{2} ( a^{-1} )
x_{1}$ \ \ \ \ \ \ \ \ \ \ \ \ \ \ (by relation \eqref{rel:5})

meanwhile,

$( a \cdot b ) \circ ( a \cdot c ) =a \tmop{sh} ( b ) s_{1} \tmop{sh}^{-1} ( a
) \tmop{sh} ( a ) \tmop{sh}^{2} ( c ) s_{2} \tmop{sh}^{2} ( a^{-1} ) x_{1}$

$=a \tmop{sh} ( b ) s_{1} \tmop{sh}^{2} ( c ) s_{2} \tmop{sh}^{2} ( a^{-1} )
x_{1}$

$=a \tmop{sh} ( b ) \tmop{sh}^{2} ( c ) s_{1} s_{2} \tmop{sh}^{2} ( a^{-1} )
x_{1}$

$=a \cdot ( b \circ c )$

Thus $\circ$ is associative and $(i)$, $(ii)$ hold.
Now, we denote 
\begin{equation}\label{eq:def_matc-B}
    \mathcal{B}= \{ a_{1} \circ   \ldots \circ a_{n} | n\in \Z_{+}, a_{i} \in
B_{\infty} \nobracket \}\subset R 
\end{equation}
and denote
\begin{equation}\label{eq:def-F_e}
    \mathcal{F}_{e}=\{ a_{1} \circ a_{2} \circ \ldots \circ a_{n} | n\in \Z_{+}, a_{i} \in
\mathcal{A}_{e} \nobracket \}\subset \mathcal{B}. 
\end{equation}
Recall that $\mathcal{A}_{e}\subset B_{\infty}$ is
the LD system generated by the unit element $e$ under composition $\cdot$ given by \eqref{eq:composition-cdot}. 
From \eqref{eq:def-compo-circ} and associativity of $\circ$, we get
\begin{equation}\label{eq:def-compo-circ-2}
    a_{1}\circ a_{2} \circ \ldots \circ a_{n}= \tmop{sh}(\vec{a})x_{1}^{n-1}
\end{equation}
for any $a_{1},\ldots,a_{n}\in R$ and $\vec{a}=(a_{1},\ldots,a_{n})$. Here, $\tmop{sh}$ is defined by equation \eqref{eq:def-sh}.

We now extend the composition $\cdot$ to $\mathcal{B}$ by equation
\eqref{eq:def-composition} in Section \ref{sec:env}. 
We need to prove that the composition $\cdot$ is well
defined on $\mathcal{B}$ and $(\mathcal{F}_{e},\cdot,\circ)$ is a free monogenerated LD monoid. The key is the following lemma whose proof is left to Section \ref{sec:order}.

\begin{lemma}\label{lem:in-B_n}
  For an element $a$ in $B_{\infty} \subseteq R$, $ax_{1}^{n-1} =x_{1}^{n-1}$
  if and only if $a \in B_{n}$.

  \begin{proof}
    We leave the proof to Section \ref{sec:order}.
  \end{proof}
\end{lemma}

\begin{theorem}\label{thm:LD-on-B}
  The composition $\cdot$ is well defined on $\mathcal{B}$ by the equation \eqref{eq:def-composition}:
\begin{multline*}
    ( a_{1} \circ a_{2} \circ \ldots \circ a_{n} ) \cdot ( b_{1} \circ b_{2}
  \circ   \ldots \circ b_{m} ) =\\ ( a_{1} \cdot ( a_{2} \cdot ( \ldots \cdot (
  a_{n} \cdot b_{1} ) ) ) \circ ( a_{1} \cdot ( a_{2} \cdot \ldots \cdot (
  a_{n} \cdot b_{2} ) ) ) \circ \ldots \circ ( a_{1} \cdot ( a_{2} \cdot
  \ldots \cdot ( a_{n} \cdot b_{m} ) ) ).
\end{multline*}
\end{theorem}
\begin{proof}
    To prove $\cdot$ is well-defined,
    suppose $a_{1} \circ a_{2} \circ   \ldots \circ a_{n}=a'_{1} \circ a'_{2} \circ   \ldots \circ a'_{n'}$ and $b_{1} \circ b_{2} \circ   \ldots b_{m}=b'_{1}\circ b'_{2} \circ   \ldots b'_{m'}$. Our goal is to prove $$
    (a_{1} \circ a_{2} \circ   \ldots \circ a_{n}) \cdot (b_{1} \circ b_{2} \circ   \ldots b_{m})=
    (a'_{1} \circ a'_{2} \circ   \ldots \circ a'_{n'}) \cdot (b'_{1}\circ b'_{2} \circ   \ldots b'_{m'}).
    $$
    By Theorem \ref{thm:B_inf-embed}, the total numbers of $x_{i}$'s are the same on both sides of these two equations. Thus by \eqref{eq:def-compo-circ-2}, we have $n=n'$ and $m=m'$. Hence, we only need to prove
    \begin{align}\label{eq:well-def-F_e}
    \begin{split}
        &( a_{1} \cdot ( a_{2} \cdot ( \ldots \cdot (
  a_{n} \cdot b_{1} ) ) ) \circ ( a_{1} \cdot ( a_{2} \cdot \ldots \cdot (
  a_{n} \cdot b_{2} ) ) ) \circ \ldots \circ ( a_{1} \cdot ( a_{2} \cdot
  \ldots \cdot ( a_{n} \cdot b_{m} ) ) )\\=
  &( a'_{1} \cdot ( a'_{2} \cdot ( \ldots \cdot (
  a'_{n} \cdot b'_{1} ) ) ) \circ ( a'_{1} \cdot ( a'_{2} \cdot \ldots \cdot (
  a'_{n} \cdot b'_{2} ) ) ) \circ \ldots \circ ( a'_{1} \cdot ( a'_{2} \cdot
  \ldots \cdot ( a'_{n} \cdot b'_{m} ) ) ).
    \end{split}
    \end{align}
    Firstly, by \eqref{eq:def-compo-circ-2}, we have $\tmop{sh}(\vec{a})x_{1}^{n-1}=\tmop{sh}(\vec{a'})x_{1}^{n-1}$ and $\tmop{sh}(\vec{b})x_{1}^{m-1}=\tmop{sh}(\vec{b'})x_{1}^{m-1}$. Here, we put $\vec{a}=( a_{1} , a_{2} , \ldots , a_{n} ), \vec{a'}=( a'_{1} , a'_{2} , \ldots , a'_{n} ), \vec{b}=( b_{1} , b_{2} ,   \ldots , b_{m} )$ and $\vec{b'}=( b'_{1} , b'_{2} ,   \ldots , b'_{m} )$.
    By Lemma \ref{lem:in-B_n}, this is equivalent to
    $\tmop{sh} ( \vec{a'} )^{-1} \tmop{sh} ( \vec{a} ) =g_{1} \in B_{n}$ and $\tmop{sh} ( \vec{b'} )^{-1} \tmop{sh} ( \vec{b} ) =g_{2} \in B_{m}$. 
    
    For any
    $c \in B_{\infty} \nocomma$, we have
    \begin{align}
        \begin{split}\label{eq:one-side}
            &a_{1} \cdot ( a_{2} \cdot ( \ldots \cdot ( a_{n} \cdot c ) ) ) \\=
            &\tmop{sh} ( \vec{a} ) \tmop{sh}^{n} ( c ) s_{n} \ldots s_{1} \tmop{sh} (
    \tmop{sh} ( \vec{a} ) )^{-1} \\=
    &\tmop{sh} ( \vec{a'} ) g_{1} \tmop{sh}^{n} ( c )
    s_{n} \ldots s_{1} \tmop{sh} ( g_{1}^{-1} ) \tmop{sh} ( \tmop{sh} ( \vec{a'} )
    )^{-1}_{}.
        \end{split}
    \end{align}
    Since $g_{1} \in B_{n}$, we have $g_{1} \tmop{sh}^{n} ( c ) = \tmop{sh}^{n} ( c )
    g_{1}$ and $s_{n} \ldots s_{1} \tmop{sh} ( g_{1}^{-1} )=g_{1}^{-1} s_{n} \ldots s_{1}$.
    Thus \eqref{eq:one-side} implies $a_{1} \cdot ( a_{2} \cdot ( \ldots \cdot ( a_{n} \cdot c ) ) )
    =a'_{1} \cdot ( a'_{2} \cdot ( \ldots \cdot ( a'_{n} \cdot c ) ) )$
    and we have 
    \begin{align}
    \begin{split}
        &( a_{1} \cdot ( a_{2} \cdot ( \ldots \cdot (
  a_{n} \cdot b_{1} ) ) ) \circ ( a_{1} \cdot ( a_{2} \cdot \ldots \cdot (
  a_{n} \cdot b_{2} ) ) ) \circ \ldots \circ ( a_{1} \cdot ( a_{2} \cdot
  \ldots \cdot ( a_{n} \cdot b_{m} ) ) )\\=
  &( a'_{1} \cdot ( a'_{2} \cdot ( \ldots \cdot (
  a'_{n} \cdot b_{1} ) ) ) \circ ( a'_{1} \cdot ( a'_{2} \cdot \ldots \cdot (
  a'_{n} \cdot b_{2} ) ) ) \circ \ldots \circ ( a'_{1} \cdot ( a'_{2} \cdot
  \ldots \cdot ( a'_{n} \cdot b_{m} ) ) ).
    \end{split}
    \end{align}
    Now, by direct calculation,
    \begin{align}
    \begin{split}
  &( a'_{1} \cdot ( a'_{2} \cdot ( \ldots \cdot (
  a'_{n} \cdot b_{1} ) ) ) \circ ( a'_{1} \cdot ( a'_{2} \cdot \ldots \cdot (
  a'_{n} \cdot b_{2} ) ) ) \circ \ldots \circ ( a'_{1} \cdot ( a'_{2} \cdot
  \ldots \cdot ( a'_{n} \cdot b_{m} ) ) )\\=
  &\tmop{sh}(\vec{a'}) \tmop{sh}^{n} ( \tmop{sh}(\vec{b})x_{1}^{m-1} ) s_{n} \ldots s_{1}
  \tmop{sh} ( \tmop{sh}(\vec{a'})^{-1} )\\=
  &\tmop{sh}(\vec{a'}) \tmop{sh}^{n} ( \tmop{sh}(\vec{b'})g_{2}x_{1}^{m-1} ) s_{n} \ldots s_{1}
  \tmop{sh} ( \tmop{sh}(\vec{a'})^{-1} )\\=
  &\tmop{sh}(\vec{a'}) \tmop{sh}^{n} ( \tmop{sh}(\vec{b'})x_{1}^{m-1} ) s_{n} \ldots s_{1}
  \tmop{sh} ( \tmop{sh}(\vec{a'})^{-1} )\\=
  &( a'_{1} \cdot ( a'_{2} \cdot ( \ldots \cdot (
  a'_{n} \cdot b'_{1} ) ) ) \circ ( a'_{1} \cdot ( a'_{2} \cdot \ldots \cdot (
  a'_{n} \cdot b'_{2} ) ) ) \circ \ldots \circ ( a'_{1} \cdot ( a'_{2} \cdot
  \ldots \cdot ( a'_{n} \cdot b'_{m} ) ) )
    \end{split}
    \end{align}
    where the third equation is due to $g_{2}x_{1}^{m-1}=x_{1}^{m-1}$ since $g_{2}\in B_{m}$ and Lemma \ref{lem:in-B_n}.
    Our theorem follows.
  \end{proof}
Note that, $\mathcal{B}= \{ a_{1} \circ   \ldots \circ a_{n} |n\in \Z_{+}, a_{i} \in
B_{\infty} \nobracket \}$ contains exactly all elements of form
$y=ax_{1}^{m}$ where $a$ is a braid. Thus Theorem \ref{thm:LD-on-B} says there is an $\tmop{LD}$
monoid structure on $\mathcal{B}$. The explicit calculation is given by:

\begin{proposition}\label{prop:explicit}
  Let $a \in B_{\infty} ,c \in \mathcal{B}$. Then we have

  (1) $( ax_{1}^{n-1} ) \cdot c=a \tmop{sh}^{n} ( c ) s_{n} \ldots s_{1}
  \tmop{sh} ( a^{-1} )$

  (2) $( ax_{1}^{n-1} ) \circ c=a \tmop{sh}^{n} ( c ) x_{1}^{n}$
\end{proposition}
\begin{proof}
  Write $a=a\circ e\circ \ldots \circ e$ where the formula contains $n-1$ $e$'s. The proposition follows from \eqref{eq:def-composition}.
\end{proof}
\begin{remark}\label{rem:rel-EB}
In \cite{dehornoy1998transfinite}, the extended braid monoid $EB_{\infty}$ is constructed where the $\circ$ operation is defined in \cite[Proposition 1.3]{dehornoy1998transfinite} and $\cdot$ operation is defined in \cite[Proposition 3.1]{dehornoy1998transfinite}. Adapting the notations in \cite{dehornoy1998transfinite}, let $EB_{\infty}^{\geq 1}=\{s \tau_{p,\infty}\in EB_{\infty}: p\geq 1\}$. Then $EB_{\infty}^{\geq 1}$ is an LD monoid. Define the map $\Theta:\mathcal{B}\rightarrow EB_{\infty}^{\geq 1}$ by $\Theta(s x_{1}^{m})= s \tau_{m+1.\infty}$. Then by direct calculation and Lemma \ref{lem:in-B_n}, $\Theta$ is an LD monoid isomorphism.
\end{remark}
\begin{corollary}\label{cor:a_cdot_b_in_R}
  For any $a,b\in \mathcal{B}$, we have $a \tmop{sh}(b) s_{1}= a\cdot b \tmop{sh}(a)$.
\end{corollary}
\begin{proof}
  By repeatedly using \eqref{eq:def-compo-circ}, we can write $a=a' x_{1}^{n-1}$ where $a'\in B_{\infty}$ and $n\in \Z_{+}$. By Proposition \ref{prop:explicit}, 
  \begin{align*}
      \begin{split}
            a\cdot b \tmop{sh}(a)=& a' \tmop{sh}^{n} ( b ) s_{n} \ldots s_{1}
  \tmop{sh} ( (a')^{-1} ) \tmop{sh}(a') x_{2}^{n-1}\\
  =& a' \tmop{sh}^{n} ( b ) s_{n} \ldots s_{1} x_{2}^{n-1}\\
  =& a' \tmop{sh}^{n} ( b ) x_{1}^{n-1} s_{1}\\
  =& a' x_{1}^{n-1} \tmop{sh}(b) s_{1}\\
  =& a \tmop{sh}(b) s_{1}.
      \end{split}
  \end{align*}
Here, we used $s_{n} \ldots s_{1} x_{2}^{n-1}=x_{1}^{n-1} s_{1}$ which is by repeatedly using rules $(3)$, $(5)$ in Proposition \ref{prop:relationofS}.
\end{proof}

\begin{theorem}\label{thm:B-is-mono}
  The $\tmop{LD}$ monoid $\mathcal{F}_{e}\subset \mathcal{B}$ defined by \eqref{eq:def-F_e}
  is the free monogenerated LD monoid with generator being the identity $e\in R$.

\end{theorem}

\begin{proof}
   Firstly, suppose $a_{1} \circ a_{2} \circ   \ldots \circ a_{m} =b_{1} \circ
    b_{2} \circ   \ldots \circ b_{n}$ with $a_{i},b_{j}\in \mathcal{A}_{e}$.

   Because $a_{1} \circ a_{2} \circ   \ldots \circ a_{m} =a_{1} \tmop{sh} ( a_{2} )
    \ldots \tmop{sh}^{m-1} ( a_{m} ) x_{1}^{m-1}$, by Theorem \ref{thm:B_inf-embed}, we have $m=n$.

    By Proposition \ref{prop:equal-in-F}, we only need to prove that
    there are elements $g_{1} ,g_{2}$ in $B_{m}^{+}$ such that $g_{1} ( a_{1}
    ,a_{2} ,  \ldots ,a_{m} ) =g_{2} ( b_{1} ,b_{2} ,  \ldots ,b_{m} )$. Recall that
    the action of positive braids is given by \eqref{eq:action-of-braid}.

    If $a_{1} \circ
    a_{2} \circ   \ldots \circ a_{m} =b_{1} \circ b_{2} \circ   \ldots \circ
    b_{m}$, then we have $$a_{1} \tmop{sh} ( a_{2} ) \ldots \tmop{sh}^{m-1} (
    a_{m} ) x_{1}^{m-1} =b_{1} \tmop{sh} ( b_{2} ) \ldots \tmop{sh}^{m-1} (
    b_{m} ) x_{1}^{m-1}.$$ Thus $\tmop{sh} ( \vec{b} )^{-1} \tmop{sh} ( \vec{a}
    ) x_{1}^{m-1} =x_{1}^{m-1}$.

    By Lemma \ref{lem:in-B_n}, $\tmop{sh} ( \vec{b} )^{-1} \tmop{sh} ( \vec{a} ) =g$ for
    some $g \in B_{m}$. Now we write $g=g_{2} g_{1}^{-1}$ where $g_{1} ,g_{2}
    \in B_{m}^{+}$. Then we have $\tmop{sh} ( \vec{a} ) g_{1} = \tmop{sh} (
    \vec{b} ) g_{2}$. Thus by Lemma \ref{lem:uniqueness-of-sh}, we have $\tmop{sh} ( g_{1} ( \vec{a} ) ) =
    \tmop{sh} ( g_{2} ( \vec{b} ) )$ and $g_{1} ( \vec{a} ) =g_{2} (
    \vec{b} )$.
  \end{proof}
\section{The order $<_{L}$ on $R$}\label{sec:order}

There is a linear ordering $\vartriangleright$ on the free group $F_{n}=\left<e_{1},\ldots,e_{n}\right>$ (and $F_{\infty}=\left<e_{1},e_{2},\ldots\right>$) that comes from hyperbolic
geometry of punctured disk. This ordering is not left-multiplicative invariant,
and is formally stated below. For geometric backgrounds which we do not use, we
refer to \cite{funk2001hurwitz},\cite{fenn1999ordering},\cite{short1999orderings},\cite{dehornoy2012braids} and especially \cite[Section 3.2 of Chapter IX]{dehornoy2008ordering}.

\begin{definition}[Definition 3.2. of Chapter IX in \cite{dehornoy2008ordering}]\label{def:order-on-free-group} 
Given $n\in \{1,2,\ldots,\infty\}$,
  there is a linear ordering on the free group $F_{n}$ defined as follows. \\
  Given two distinct elements $w$, $u$ in $F_{n}$, write them in the minimal word form, say, $w=e_{i_{1}}^{k_{1}} \ldots e_{i_{l}}^{k_{l}}$ and
  $u=e_{j_{1}}^{q_{1}} \ldots e_{j_{p}}^{q_{l}}$ where $k_{h} ,q_{h} = \pm 1$.
  Suppose $m$ satisfies $e_{i_{r}}^{k_{r}} =e_{j_{r}}^{q_{r}}$ when $r<m$ and
  $e_{i_{m}}^{k_{m}} \neq e_{j_{m}}^{q_{r}}$. Then the following holds.
\begin{itemize}
    \item If $m=1$, $w \vartriangleleft u$ if and only if $k_{1} i_{1} >q_{1} j_{1}$. 
    \item For $m>1$ and $m \leqslant \min ( l,p )$:

  (1) If $k_{m-1} =1$, $w \vartriangleleft u$ if only if one of the
  following holds:\\
  a. $k_{m} ,q_{m} =1$ and $i_{m} >j_{m}$\\ b. $k_{m} =1,q_{m} =-1$ and $j_{m}
  <j_{m-1}$\\ c. $k_{m} =-1,q_{m} =1$ and $i_{m} >i_{m-1}$ \ \\d. $k_{m} ,q_{m}
  =-1$ and $( ( i_{m} >i_{m-1}   \tmop{and}  i_{m} >j_{m}   )  \; \tmop{or}\;   (
  i_{m} <j_{m} <i_{m-1} )   )$. \

  (2) If $k_{m-1} =-1$, $w \vartriangleleft u$ if and only if $\iota ( u )
  \vartriangleleft \iota ( w )$ where $\iota$ is the automorphism of $F_{n}$
  mapping $e_{i}$ to $e_{i}^{-1}$ for each $i\geq 1$.
  \item If $m-1=p<l$, $w \vartriangleleft u$ if and only if $( \nobracket k_{m-1}
  ,k_{m} =-1$ and $i_{m} >i_{m-1} \nobracket )  $ or
  $( k_{m-1} =1  \text{ and }   ( k_{m} =-1  \tmop{or}  i_{m} <i_{m-1} ) )$.
\end{itemize}
\end{definition}

\begin{proposition}[Section 3 of Chapter IX in \cite{dehornoy2008ordering}]\label{prop:action-braid-on-free}
  Let $n \in \{ 1,2, \ldots , \infty \}$. There is a group action of
  $B_{n}$ as outer automorphisms on $F_{n}$ given by:
  \[ s_{i} ( e_{i} ) =e_{i+1} e_{i}^{-1} e_{i-1} \text{ , } s^{-1}_{i} ( e_{i} ) =e_{i-1} e_{i}^{-1} e_{i+1} \]
  for $i=1,2, \ldots $ where $e_{0} =e$,
  \[ s_{i} ( e_{j} ) =e_{j} \text{ , } s^{-1}_{i} ( e_{j} ) =e_{j} \]
  for $i \neq j$.
\end{proposition}

\begin{proposition}[Proposition 3.4 of Chapter IX in \cite{dehornoy2008ordering}]\label{prop:order-preserving}
  The action given by Proposition 6.2 is $\vartriangleleft -
  \tmop{preserving}$.
\end{proposition}

\begin{lemma}\label{lem:derhornoy-order-rep}
  For any two elements $a,b \in
B_{\infty}$, $a<_{D} b$ if and only if there is an $n>0$, such that $a ( e_{n}
) \vartriangleleft b ( e_{n} )$ and $a ( e_{i} ) =b ( e_{i} )   \forall i<n$ (recall that $<_{D}$ is defined in Theorem \ref{thm:linear-Dehornoy}).
\end{lemma}
\begin{proof}
    Since $<_{D}$ is linear on $B_{\infty}$, we only need to prove that, $a<_{D} b$ implies $a ( e_{n} ) \vartriangleleft b ( e_{n} )$ and $a ( e_{i} ) =b ( e_{i} )   (\forall i<n)$ for some $n\in \Z_{+}$. 
    To see this, suppose $a<_{D} b$. By Theorem \ref{thm:linear-Dehornoy}, there is $n\in \Z_{+}$ and $c_{i}\in \langle s_{n+1},s_{n+2},\ldots \rangle (1\leq i\leq m)$ such that $$a^{-1} b= c_{1}s_{n}c_{2}s_{n}\ldots c_{m-1} s_{n} c_{m}.$$
    Consider the action of $B_{\infty}$ on $F_{\infty}$ defined in Proposition \ref{prop:action-braid-on-free}. Then under this action, $c_{i} (i\leq m)$ and $s_{n}$ all fix $e_{j} (j<n)$.
    Thus we have $a^{-1}b (e_{j})=e_{j}$ for $j<n$. Hence it suffices to prove $a^{-1}b(e_{n}) \vartriangleright e_{n}$. By Definition \ref{def:order-on-free-group}, we have $s_{n}(e_{n})=e_{n+1}e_{n}^{-1}e_{n-1} \vartriangleright e_{n}$. Together with Proposition \ref{prop:order-preserving} and $c_{i}(e_{n})=e_{n} (i\leq m)$, we have
    \begin{align}
        \begin{split}
            a^{-1} b (e_{n})=& c_{1}s_{n}c_{2}s_{n}\ldots c_{m-1} s_{n} c_{m} (e_{n})\\
            =& c_{1}s_{n}c_{2}s_{n}\ldots c_{m-1} s_{n} (e_{n})\\
            \vartriangleright& c_{1}s_{n}c_{2}s_{n}\ldots c_{m-1} (e_{n})\\
            \trianglerighteqslant&\cdots \\
            \trianglerighteqslant& c_{1} (e_{n})\\
            =& e_{n}.
        \end{split}
    \end{align}
    By Proposition \ref{prop:order-preserving} again, we have $a(e_{n})\vartriangleleft b(e_{n})$ and
    the lemma follows.
\end{proof}



Lemma \ref{lem:derhornoy-order-rep} gives an equivalent description of Dehornoy order by using representation of $B_{\infty}$ as automorphisms of $F_{\infty}$ defined in Proposition \ref{prop:action-braid-on-free}.
Now, we will represent monoid $R$ as injective morphisms from $F_{\infty}$ to itself (which extends the action of $B_{\infty}$ on $F_{\infty}$) and use this representation to extend the Dehornoy order from $B_{\infty}$ to $R$.

Let us denote by $Hom_{Inj}(F_{\infty},F_{\infty})$ the monoid of injective order-preserving group morphisms from $( F_{\infty} ,
  \vartriangleleft )$ to itself where the monoid composition on $Hom_{Inj}(F_{\infty},F_{\infty})$ is given by the composition of functions.
\begin{proposition}\label{prop:action-R-on-F_infty}
  There is a monoid morphism from $R$ to $Hom_{Inj}(F_{\infty},F_{\infty})$ which induces an order-preserving action of $R$ on $F_{\infty}$ such that:
  \[ s_{i} ( e_{i} ) =e_{i+1} e_{i}^{-1} e_{i-1} \text{ , } s^{-1}_{i} ( e_{i} ) =e_{i-1} e_{i}^{-1} e_{i+1} \]
  for $i=1,2, \ldots $ where $e_{0} =e$,
  \[ s_{i} ( e_{j} ) =e_{j} \text{ , } s^{-1}_{i} ( e_{j} ) =e_{j} \]
  for $i \neq j$.
  Also, $x_{i} ( e_{j} ) =e_{j+1}$ for $j \geqslant i$, \ $x_{i} ( e_{j} )
  =e_{j}$ for $j<i$.
\end{proposition}
\begin{proof}
    We need to prove the relations $( 1 ) \sim ( 7 )$ in Proposition \ref{prop:relationofS} hold after substituting $\tilde{s_{i}},\tilde{x_{i}}$ by corresponding elements in $Hom_{Inj}(F_{\infty},F_{\infty})$ that are defined in the conditions of Proposition \ref{prop:action-R-on-F_infty}. Note that $( 6 ) ,
    ( 7 )$ hold because of Proposition \ref{prop:action-braid-on-free}. We will check the relation $( 1 )
    \sim ( 3 )$ hold and leave $(4)$, $(5)$ to the reader.

    \subparagraph{(1)$s_{i} x_{i} =x_{i}$:}

    $s_{i} x_{i} ( e_{j} ) =s_{i} ( e_{j+1} ) =e_{j+1} =x_{i} ( e_{j} )$ for
    $j \geqslant i$ and $s_{i} x_{i} ( e_{j} ) =s_{i} ( e_{j} ) = e_{j} =s_{i}
    ( e_{j} )$ for $j<i$.

    \subparagraph{(2)$x_{i+1} x_{i} =x_{i}^{2}$:}

    $x_{i+1} x_{i} ( e_{j} ) =e_{j+2} =x_{i}^{2} ( e_{j} )$ for $j \geqslant
    i$ and $x_{i+1} x_{i} ( e_{j} ) =e_{j} =x_{i}^{2} ( e_{j} )$ for $j<i$.

    \subparagraph{(3)$ x_{i+1} s_{i} =s_{i} s_{i+1} x_{i}$:}

    $x_{i+1} s_{i} ( e_{i} ) =x_{i+1} ( e_{i+1} e_{i}^{-1} e_{i-1} ) =e_{i+2}
    e_{i}^{-1} e_{i-1}$ and
    \begin{align*}
            s_{i} s_{i+1} x_{i} ( e_{i} ) =&s_{i} s_{i+1} (
    e_{i+1} ) =s_{i} ( e_{i+2} e_{i+1}^{-1} e_{i} )\\
    =&e_{i+2} e_{i+1}^{-1} e_{i+1} e_{i}^{-1} e_{i-1} =e_{i+2} e_{i}^{-1}
    e_{i-1}.
    \end{align*}
    For $j>i$, $x_{i+1} s_{i} ( e_{j} ) =x_{i+1} ( e_{j} ) =e_{j+1} =s_{i}
    s_{i+1} x_{i} ( e_{j} )$ .

    As for the order preserving property, this is true for $s_{i}$ by Proposition \ref{prop:order-preserving}, and is true for $x_{i}$ by directly checking Definition \ref{def:order-on-free-group}.
  \end{proof}

Now we prove Lemma \ref{lem:in-B_n} in Section \ref{sec:monoidR}:
\begin{proof}
{\tmem{Suppose}}{\tmem{ $ax_{1}^{n-1} =x_{1}^{n-1}$ in $R$ and assume $a \nin
B_{n}$. By Theorem \ref{thm:word-form-in-B_inf}, there is $k\in \Z_{+}$ such that we can write $a$ in a word form
of $s_{1}^{\pm} , \ldots ,s_{k}^{\pm}$ where one of $s_{k}$ and $s_{k}^{-1}$
appears. Now, if $k<n$ then $a \in B_{n}$ thus we must have $k \geqslant n$.
}}

{\tmem{Consider the element $e_{k-n+1}$ in $F_{\infty}$. By Proposition \ref{prop:action-R-on-F_infty}, we have $$x_{1}^{n-1}
( e_{k-n+1} ) =e_{k}.$$ One the other hand, $ax_{1}^{n-1} ( e_{k-n+1} ) =a (
e_{k} ) \neq e_{k}$ by using the argument in the proof of Lemma \ref{lem:derhornoy-order-rep}. Arrive at contradiction.}}
\end{proof}
\begin{definition}
  Define the partial order $<_{L}$ on $R$ as follows: for any two elements $a,b \in
  R$, $a<_{L} b$ if and only if there is an $n \in \Z_{+}$, such that $a ( e_{n} ) \vartriangleleft b (
  e_{n} )$ and $a ( e_{i} ) =b ( e_{i} )$ $\forall i<n$,
\end{definition}

\begin{proposition}\label{prop:linearity_<_L}
  The order $<_{L}$ is left-invariant on $R$. Restricting
  on $B_{\infty}$, $<_{L}$ is the Dehornoy order $<_{D}$.
\end{proposition}
\begin{proof}
    Since $R$ acts on $F_{\infty}$ as order preserving morphisms, $<_{L}$ is
    left-invariant. When restricting on $B_{\infty}$, $<_{L}$ is the Dehornoy
    order by Lemma \ref{lem:derhornoy-order-rep}.
  \end{proof}
Now we prove that $<_{L}$ is a linear order. Then we deduce that the morphism from $R$ to $Hom_{Inj}(F_{\infty},F_{\infty})$ defined in Proposition \ref{prop:action-R-on-F_infty} is injective.

In Section \ref{sec:monoidR}, we defined in \eqref{eq:def-compo-circ} a composition $\circ$ on $R$. Together with composition $\cdot$ (defined in \eqref{eq:def-composition}), we endowed on
$\mathcal{B}$ (defined in \eqref{eq:def_matc-B}) an LD monoid structure. By Theorem \ref{thm:B-is-mono}, $\mathcal{F}_{e}\subset \mathcal{B}$ is a free
monogenerated $\tmop{LD}$ monoid. Given $a,b \in \mathcal{F}_{e}$, we write
$a<b$ if $a<b$ in the Laver order defined in Theorem \ref{thm:linear-Laver}.
\begin{lemma}\label{lem:coincide-on-F_e}
  For any two $a,b\in \mathcal{F}_{e}$, $a(e_{1})\vartriangleleft b(e_{1})$ if and only if $a<b$. 
\end{lemma}
\begin{proof}
    We only need to prove $a \cdot b (e_{1}) \vartriangleright a (e_{1})$ and $a\circ b (e_{1})\vartriangleright a (e_{1})$ for any $a,b\in \mathcal{F}_{e}$, then our conclusion follows by Theorem \ref{thm:linear-Laver}. We first claim that $a(e_{1}) \vartriangleright e_{1}$ for any $a\in \mathcal{F}_{e}$ with $a\not= e$. To see this, write 
    $a=a_{1} \circ a_{2} \circ \ldots \circ a_{n}$ with $a_{i}\in \mathcal{A}_{e}$. If $n=1$, then $a\in \mathcal{A}_{e}$ and our claim follows from the proof of Lemma \ref{lem:derhornoy-order-rep}. Otherwise, 
    \begin{align*}
        a (e_{1})
        =& a_{1} \tmop{sh}(a_{2})\ldots \tmop{sh}^{n-1} (a_{n}) x_{1}^{n-1} (e_{1})\\
        =& a_{1} \tmop{sh}(a_{2})\ldots \tmop{sh}^{n-1} (a_{n}) (e_{n})\\
        \vartriangleright& a_{1} \tmop{sh}(a_{2})\ldots \tmop{sh}^{n-1} (a_{n}) (e_{1})\\
        =& a_{1} (e_{1})\\
        \trianglerighteqslant& e_{1}.
    \end{align*}
    Here, the first inequality is due to $e_{n}\vartriangleright e_{1}$ when $n>1$ and the fact that the action of $R$ preserves the order $\vartriangleright$. Thus our claim holds.
    
    Now, suppose $a,b \in \mathcal{F}_{e}$ with $a=a_{1} \circ a_{2} \circ
    \ldots \circ a_{n}$, and $b=b_{1} \circ b_{2} \circ   \ldots \circ b_{m}$
    where $a_{i} ,b_{j} \in \mathcal{A}_{e}$. We let $\vec{a}=(a_{1},a_{2},\ldots,a_{n})$.
    By Proposition \ref{prop:explicit}, we have
    \begin{equation}
        a \cdot b
    = \tmop{sh} ( \vec{a} ) \tmop{sh}^{n} ( b ) s_{n} \ldots s_{1} \tmop{sh} (
    \tmop{sh} ( \vec{a} )^{-1} )
    \end{equation} and
    \begin{align*}
        a \cdot b ( e_{1} ) =& \tmop{sh} ( \vec{a} ) \tmop{sh}^{n} ( b ) s_{n}
    \ldots s_{1} ( e_{1} ) \\
    \vartriangleright& \tmop{sh} ( \vec{a} )
    \tmop{sh}^{n} ( b ) (e_{n})\\
    =& \tmop{sh} ( \vec{a} )
    (e_{n})\\
    =& a_{1} \circ a_{2} \circ   \ldots \circ a_{n} (
    e_{1} )\\
    =& a (e_{1}).
    \end{align*}
    Here, the first inequality is because $s_{n} \cdot \ldots \cdot s_{1} ( e_{1} )=e_{n+1} e_{n}^{-1} 
    \vartriangleright e_{n}$ by checking Definition \ref{def:order-on-free-group}. Thus we have $a \cdot b (e_{1})\vartriangleright a(e_{1})$.
    
    It remains to prove $a \circ b (e_{1})\vartriangleright a(e_{1})$. By Proposition \ref{prop:explicit}, we have 
    \begin{align*}
        a \circ b (e_{1}) =& \tmop{sh}(\vec{a}) \tmop{sh}^{n}(b) x_{1}^{n} (e_{1})\\
        =& \tmop{sh}(\vec{a}) \tmop{sh}^{n} (b) (e_{n+1})\\
        \trianglerighteqslant & \tmop{sh}(\vec{a}) (e_{n+1})\\
        \vartriangleright & \tmop{sh}(\vec{a}) (e_{n})\\
        = & \tmop{sh}(\vec{a}) x_{1}^{n-1} (e_{1})\\
        = & a(e_{1}).
    \end{align*}
    Here, we used $b(e_{1})\trianglerighteqslant e_{1}$ by our claim above and thus $\tmop{sh}^{n} (b) (e_{n+1})\trianglerighteqslant e_{n+1}$. We also used $e_{n+1}\vartriangleright e_{n}$ for each $n\in \Z_{+}$ by Definition \ref{def:order-on-free-group}.
\end{proof}
\begin{definition}\label{def:order-on-seq}
Given $n\in \Z_{+}$, we define $\mathcal{F}^{(n)}=\{(a_{1},\ldots,a_{n}):a_{i}\in \mathcal{F}_{e}\}$ to be the set of sequences in $\mathcal{F}_{e}$ with length $n$. For any $n,m \in \Z_{+}$ and $\vec{a}\in \mathcal{F}^{(n)}$ and $\vec{b}\in \mathcal{F}^{(m)}$, we define $\vec{a} <_{s} \vec{b}$ if and only if one of following holds:
\begin{itemize}
    \item $n<m$ and $\vec{a}$ is a prefix of $\vec{b}$
    \item there exists $k\leq \min\{m,n\}$
such that $a_{i} =b_{i}$ for $i<k$ and $a_{k} <b_{k}$ in $\mathcal{F}_{e}$.
\end{itemize}
Note that $<_{s}$ is a linear order on $\bigcup_{n\geq 1} \mathcal{F}^{(n)}$ since $<$ is linear on $\mathcal{F}_{e}$.

For any finite sequence $\vec{a}$ of elements in $\mathcal{F}_{e}$, we also denote $$\tmop{sh}(\vec{a})=a_{1} \tmop{sh}(a_{2}) \ldots \tmop{sh}^{n-1}(a_{n})$$ which is in $R$.
\end{definition}

\begin{theorem}\label{thm:uniqueness-F}
  Suppose $\vec{a}=(a_{1},\ldots,a_{n})$ and $\vec{b}=(b_{1},\ldots,b_{m})$ are two sequences of
  $\mathcal{F}_{e}$ with $a_{n} ,b_{m} \neq e$. Then $\tmop{sh}(\vec{a})<_{L} \tmop{sh}(\vec{b})$ if and only if $\vec{a}<_{s}\vec{b}$. In particular, for any $y\in R$, the representation $y= \tmop{sh} ( \vec{a} )$ with $a_{i} \in
  \mathcal{F}_{e}$ and $a_{n} \neq e$ is unique. i.e. If there are
  $b_{j} \in \mathcal{F}_{e}$ for $j \leqslant m$ and $b_{m} \neq e$ such that $y=
  \tmop{sh} ( \vec{b} )$, then $\vec{a} = \vec{b}$.
\end{theorem}
\begin{proof}
  By linearity of Laver order, $\vec{a}\not=\vec{b}$ implies $\vec{a}<_{s}\vec{b}$ or $\vec{b}<_{s}\vec{a}$. Without loss of generality, we assume $\vec{a}<_{s}\vec{b}$. Our goal is to prove $\tmop{sh}(\vec{a})<_{L} \tmop{sh}(\vec{b})$. 
  
  Case 1: $\vec{a}$ is a prefix of $\vec{b}$. In this case, let $k$ be the minimal integer with $k>n$ and $b_{k}\not=e$. Then $b_{k}>e$ and $b_{k}(e_{1})\vartriangleright e_{1}$ by Lemma \ref{lem:coincide-on-F_e} and thus $\tmop{sh}^{k-1}(b_{k})(e_{k}) \vartriangleright e_{k}$. This implies
  \begin{align*}
      \begin{split}
          \tmop{sh}(\vec{b})(e_{k}) = \tmop{sh}(\vec{a}) \tmop{sh}^{n}(b_{n+1})\ldots \tmop{sh}^{m-1}(b_{m})(e_{k})=
           \tmop{sh}(\vec{a}) \tmop{sh}^{k-1}(b_{k})(e_{k})
          \vartriangleright \tmop{sh}(\vec{a}) (e_{k}).
      \end{split}
  \end{align*}
  Meanwhile, for $j<k$, $\tmop{sh}(\vec{b})(e_{j})=\tmop{sh}(\vec{a}) (e_{j})$ by definition of $k$. Thus we have $\tmop{sh}(\vec{a})<_{L} \tmop{sh}(\vec{b})$. 
  
  Case 2: There exists $k'$ with $a_{k'}<b_{k'}$ and $a_{j}=b_{j}$ for $j<k'$. In this case, 
  $$ \tmop{sh}(\vec{a})(e_{j})= \tmop{sh}(\vec{b}) (e_{j})
  $$
  for $j<k'$. Meanwhile, we have $a_{k'}(e_{1})\vartriangleleft b_{k'}(e_{1})$ by Lemma \ref{lem:coincide-on-F_e} and thus $\tmop{sh}^{k'-1}(a_{k'})(e_{k'}) \vartriangleleft \tmop{sh}^{k'-1}(b_{k'})(e_{k'})$ and
  \begin{align*}
      \begin{split}
          \tmop{sh}(\vec{a})(e_{k'})=&a_{1} \tmop{sh}(a_{2}) \ldots \tmop{sh}^{k'-1}(a_{k'})(e_{k'})\\
          \vartriangleleft& a_{1} \tmop{sh}(a_{2}) \ldots \tmop{sh}^{k'-1}(b_{k'})(e_{k'})\\
          =&\tmop{sh}(\vec{b}) (e_{k'}).
      \end{split}
  \end{align*}
  \end{proof}
\begin{definition}
  Given $1\leq i<n$, define function $\eta_{s_{i}}:\mathcal{F}^{(n)}\rightarrow \mathcal{F}^{(n)}$ by
\begin{equation}
    \eta_{s_{i}} (a_{1},\ldots,a_{i},a_{i+1},\ldots,a_{n})=(a_{1},\ldots,a_{i}\cdot a_{i+1},a_{i},\ldots,a_{n}),
\end{equation}
we also define $\eta_{x_{i}}:\mathcal{F}^{(n)}\rightarrow \mathcal{F}^{(n-1)}$ by 
\begin{equation}
    \eta_{x_{i}} (a_{1},\ldots,a_{i-1},a_{i},a_{i+1},\ldots,a_{n})=(a_{1},\ldots,a_{i-1},a_{i}\circ a_{i+1},a_{i+2},\ldots,a_{n}).
\end{equation}
\end{definition}
We remark that, the function $\eta_{s_{i}}$ is a generalization of the action of $s_{i}$ on $\mathfrak{B}^{(n)}$. See equation \eqref{eq:action-of-braid}.
\begin{corollary}\label{cor:linear-order-R+}
  $<_{L}$ restricting on $R^{+}$ is a linear order.
\end{corollary}
\begin{proof}
    We first claim that, for any $\vec{a}\in \mathcal{F}^{(n)}$ and any $1\leq i<n$, we have
    \begin{equation}\label{eq:s_i-action}
        \tmop{sh} ( \eta_{s_{i}} ( \vec{a} ) ) = \tmop{sh} ( \vec{a} ) s_{i}
    \end{equation}
    and
    \begin{equation}\label{eq:x_i-action}
            \tmop{sh} ( \eta_{x_{i}} ( \vec{a} ) ) = \tmop{sh} ( \vec{a} ) x_{i}.
    \end{equation}

    They are calculated as follows:
\begin{align*}
    \tmop{sh} ( \vec{a} ) s_{i} =&a_{1} \tmop{sh} ( a_{2} ) \ldots \tmop{sh}^{n-1}
( a_{n} ) s_{i} \\=&a_{1} \tmop{sh} ( a_{2} ) \ldots \tmop{sh}^{i-1} ( a_{i} )
\tmop{sh}^{i} ( a_{i+1} ) s_{i} \tmop{sh}^{i+1} ( a_{i+2} ) \ldots
\tmop{sh}^{n-1} ( a_{n} )
\\=&a_{1} \tmop{sh} ( a_{2} ) \ldots \tmop{sh}^{i-1} ( a_{i} \tmop{sh} ( a_{i+1}
) s_{1} ) \tmop{sh}^{i+1} ( a_{i+2} ) \ldots \tmop{sh}^{n-1} ( a_{n} )
\\=&a_{1} \tmop{sh} ( a_{2} ) \ldots \tmop{sh}^{i-1} ( a_{i} \cdot a_{i+1} \tmop{sh}(a_{i}) ) 
\tmop{sh}^{i+1} ( a_{i+2} ) \ldots \tmop{sh}^{n-1} ( a_{n} )
\\=&a_{1} \tmop{sh} ( a_{2} ) \ldots \tmop{sh}^{i-1} ( a_{i} \cdot a_{i+1} )
\tmop{sh}^{i} ( a_{i} ) \tmop{sh}^{i+1} ( a_{i+2} ) \ldots \tmop{sh}^{n-1} (
a_{n} )
\\= &\tmop{sh} ( \eta_{s_{i}} ( \vec{a} ) )
\end{align*}
where the third equation is due to Corollary \ref{cor:a_cdot_b_in_R}.
\begin{align*}
    \tmop{sh} ( \vec{a} ) x_{i} =&a_{1} \tmop{sh} ( a_{2} ) \ldots \tmop{sh}^{n-1}
( a_{n} ) x_{i} \\=&a_{1} \tmop{sh} ( a_{2} ) \ldots \tmop{sh}^{i-1} ( a_{i} )
\tmop{sh}^{i} ( a_{i+1} ) x_{i} \tmop{sh}^{i} ( a_{i+2} ) \ldots
\tmop{sh}^{n-2} ( a_{n} )
\\=&a_{1} \tmop{sh} ( a_{2} ) \ldots \tmop{sh}^{i-1} ( a_{i} \tmop{sh} ( a_{i+1}
) x_{1} ) \tmop{sh}^{i} ( a_{i+2} ) \ldots \tmop{sh}^{n-2} ( a_{n} )
\\=&a_{1} \tmop{sh} ( a_{2} ) \ldots \tmop{sh}^{i-1} ( a_{i} \circ a_{i+1} )
\tmop{sh}^{i} ( a_{i+2} ) \ldots \tmop{sh}^{n-2} ( a_{n} ) \;\;\;\;\;\;\;\;\;\;\;\;(\text{by \eqref{eq:def-compo-circ}})
\\=& \tmop{sh} ( \eta_{x_{i}}
( \vec{a} ) ).\\
\end{align*}

    Now by definition, any element $b$ in $R^{+}$ can be written in a word form of $x_{i}'s$ and ${s_{i}}'s$. Thus we can use \eqref{eq:s_i-action} and \eqref{eq:x_i-action} repeatedly to write $b=\tmop{sh}(\vec{a})$ for some $\vec{a}\in \mathcal{F}^{(n)}$ where $n\in \Z_{+}$. More precisely, we first write $b=u_{1} u_{2}\ldots u_{t}$ where $u_{j}\in \bigcup_{i=1}^{M}\{s_{i},x_{i}:i\geq 1\}$ for each $1\leq j\leq t$ and $M$ is a sufficiently large integer. Denote $\vec{e}=(e,e,\ldots,e)\in \mathcal{F}^{(M+1)}$ where each entry is the identity $e$. Then by letting $\vec{a}=\eta_{u_{t}} \eta_{u_{t-1}} \ldots \eta_{u_{1}} (\vec{e})$, we have
    \begin{align*}
     b=& \tmop{sh}(\vec{e}) u_{1} u_{2} \ldots u_{t}\\
     =& \tmop{sh}(\eta_{u_{1}}(\vec{e})) u_{2}\ldots u_{t}\\
     =& \ldots \\
     =& \tmop{sh}(\eta_{u_{t}} \eta_{u_{t-1}} \ldots \eta_{u_{1}} (\vec{e}))\\
     =& \tmop{sh}(\vec{a}).
    \end{align*}
    If $b\not=e$, we can also assume the last entry $a_{n}\not =e$. By Theorem \ref{thm:uniqueness-F}, the vector $\vec{a}\in \mathcal{F}^{(n)}$ with $a_{n}\not=e$ and $\tmop{sh}(\vec{a})=b$ is unique. Thus the linearity of $<_{L}$ on $R^{+}$ follows from Theorem \ref{thm:uniqueness-F} and the linearity of $<_{s}$ on $\bigcup_{n\geq 1} \mathcal{F}^{(n)}$. 
  \end{proof}

\begin{theorem}\label{thm:beta-injective}
  The order $<_{L}$ on $R$ is a linear order and the morphism defined in
  Proposition \ref{prop:action-R-on-F_infty} from $R$ to $Hom_{Inj} ( F_{\infty} ,F_{\infty} )$ is
  injective.

  \begin{proof}
    In order to prove $<_{L}$ is linear, it suffices to prove that, for any $y\not=z$ in $R$, one of $y<_{L} z$ and $z<_{L} y$ holds. 
    Given $y\not=z$ in $R$,
    we can use rule $( 3 ) \sim ( 5 )$ in Proposition \ref{prop:relationofS} to move $x_{i}$'s to the tail of the word forms of $y,z$.
    Thus we can write $y=s x$ and $z=s' x'$ where $s,s' \in B_{\infty}$ and $x,x' \in X_{\infty}$. Let $s^{-1} s'=s^{-1}_{*} s_{**}$ where $s_{*},s_{**}\in B^{+}_{\infty}$. Then $y<_{L} z$ is equivalent to $s_{*} x <_{L} s_{**} x'$ and $z <_{L} y$ is equivalent to $s_{**} x <_{L} s_{*} x'$. By Corollary \ref{cor:linear-order-R+}, since $s_{*}x, s_{**}x' \in R^{+}$ and $<_{L}$ is linear on $R^{+}$, one of $s_{*}x <_{L} s_{**}x'$ and $s_{**}x'<_{L} s_{*}x$ holds. Thus $<_{L}$ is a linear order on $R$. 
    
    Now, for any $y_{1}\not = y_{2}$ in $R$, let $f_{y_{1}},f_{y_{2}} \in Hom_{Inj} ( F_{\infty} ,F_{\infty} )$ be the images of $y_{1},y_{2}$ under the morphism defined in Proposition \ref{prop:action-R-on-F_infty}.
    Assume without loss of generality that $y_{1}<_{L} y_{2}$. Then by Definition of $<_{L}$, there exists $i \geq 1$ with $f_{y_{1}} (e_{i})\vartriangleleft f_{y_{2}} (e_{i})$ and thus $f_{y_{1}}\neq f_{y_{2}}$. This proves that the morphism in Proposition \ref{prop:action-R-on-F_infty} is injective.
  \end{proof}
\end{theorem}

We denote by $\beta$ the injective morphism from $R$ to $Hom_{Inj} ( F_{\infty}
,F_{\infty} )$ defined in Proposition \ref{prop:action-R-on-F_infty}.

\section{The relations of $S$}\label{sec:relations}
\begin{figure}
    \centering
    \includegraphics{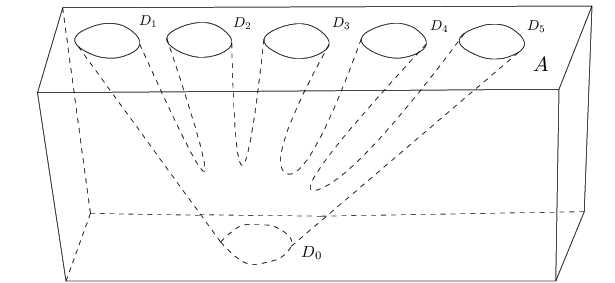}
    \caption{The picture illustrates a $5$-pant embedded in $M$.}
    \label{fig:j-pant}
\end{figure}
Recall that in Proposition \ref{prop:define-alpha}, we defined a morphism $\alpha$ from $R$ to $S$ where $S$ is the set of equivalent classes of
cobordisms between two countable collections of circles generated by $\alpha (
x_{i} )=\tilde{x_{i}}$ and $\alpha ( s_{i} )=\tilde{s_{i}}$. In this section, we prove Theorem \ref{thm:main}, namely, $\alpha:R \rightarrow S$ is an isomorphism.

The strategy is to construct a function $\phi$ from $S$ to $\tmop{Hom} (
F_{\infty} ,F_{\infty} )$ such that $\phi \circ \alpha = \beta$. Since by Theorem \ref{thm:beta-injective}, $\beta$ is injective, we imply that $\alpha$ is injective. Together with Proposition \ref{prop:define-alpha}, we conclude that $\alpha$ is bijective and thus Theorem \ref{thm:main} holds. 

We first consider the a simpler case when the cobordism does not have braid components.
\begin{definition}\label{def:j-pant}
  Let $j$ be a positive integer and cuboid $M=[0,1]^{3}$. Suppose we have $j$ input disks $D_{k} (1\leq k\leq j)$ on the upper surface $A=\{(x,y,1):0\leq x,y\leq 1\}$ and one output disk $D_{0}$ on the lower surface $\{(x,y,0):0\leq x,y\leq 1\}$.
  We call the corbodism $\rho$ embedded in $M$ with input disks $D_{k} (1\leq k\leq j)$ and output disk $D_{0}$ a $j$-pant. See Figure \ref{fig:j-pant} for an illustration. The (topological) surface $\Bar{\rho}= \rho \cup \bigcup_{0\leq k\leq j}D_{k}$ separates $M$ into two connected components $M\setminus \Bar{\rho}=M^{(1)}\cup M^{(2)}$. Let $M(\rho)$ denote the connected component which contains some boundary points of $M$ (i.e. the exterior of surface $\Bar{\rho}$).  With a little abuse of notations, we also call $M(\rho)$ the exterior of $\rho$.
\end{definition}

\begin{lemma}\label{lem:j-pant-lem}
  Consider a j-pant $\rho$ embedded in a cuboid $M$ in the way described in
  definition \ref{def:j-pant}. Let $\tilde{A} =A  \backslash
  \cup_{l=1}^{j} D_{l} \nobracket$. Then the inclusion $i: \tilde{A}
  \hookrightarrow M ( \rho )$ induces an isomorphism $i_{\ast} : \pi_{1} (
  \tilde{A} ) \rightarrow \pi_{1} ( M ( \rho ) )$.
  \end{lemma}
  \begin{proof}
    We prove by induction on $j$.

    For $j=1$, this is true because $\rho$ is a tube in $M$ and there is an
    obvious contraction from $M ( \rho )$ to $\tilde{A}$.

    Suppose this holds for $j-1$. Consider a surface $\gamma$ ``splitting''
    $M ( \rho )$ into two parts $U_{1}$ and $U_{2}$ such that $U_{1}\cap U_{2}=\gamma\cap M(\rho)$ and $U_{1}\cup U_{2}=M(\rho)$ in the way shown in Figure \ref{fig:U_1andU_2}.
\begin{figure}
    \centering
    \includegraphics{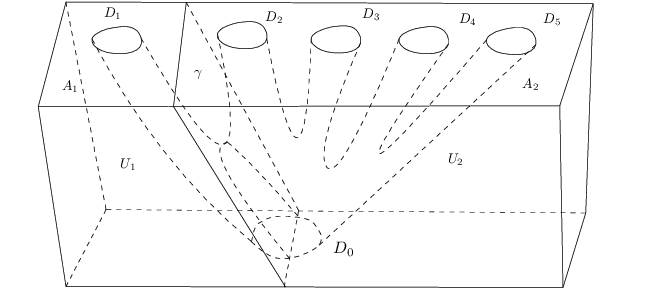}
    \caption{$\rho$ is a $j$-pant and the affine plane $\gamma$ separates $M(\rho)$ into $U_1$ and $U_2$ where $U_2$ is essentially the exterior of a $(j-1)$-pant.}
    \label{fig:U_1andU_2}
\end{figure}
    Let
    $U_{1} \cap \tilde{A} =A_{1}$, $U_{2} \cap \tilde{A} =A_{2}$ and $A_{1} \cap
    A_{2} = \tau$. Note that $U_{2}$ is essentially the exterior of a $(
    j-1 )$-pant embedded in a cuboid. Thus by induction hypothesis, $i_{2} :A_{2} \hookrightarrow
    U_{2}$ induces an isomorphism $i_{2 \ast} : \pi_{1} ( A_{2} ) \rightarrow
    \pi_{1} ( U_{2} )$. Meanwhile, the case when $j=1$ gives the
    isomorphism $i_{1 \ast} : \pi_{1} ( A_{1} ) \rightarrow \pi_{1} ( U_{1} )$
    induced by inclusion $i_{1} :A_{1} \hookrightarrow U_{1}$. Evidently,
    $i_{0 \ast} : \pi_{1} ( \tau ) \rightarrow \pi_{1} ( \gamma\cap M(\rho) )$ induced by
    inclusion $i_{0} : \tau \hookrightarrow \gamma\cap M(\rho)$ is an isomorphism between
    trivial groups. Now by using the van Kampen theorem and the naturality of
    pushout, we deduce that $i_{\ast}$ is an isomorphism. Thus the induction principle proves the lemma.
  \end{proof}

\begin{figure}
    \centering
    \includegraphics{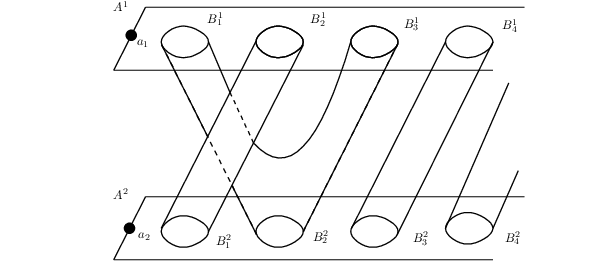}
    \caption{The picture illustrates the shrinking braid $\widetilde{s_{1}} \widetilde{x_{2}}$.}
    \label{fig:shrinking-braid}
\end{figure}

Recall notations from Section \ref{sec:shrink}.
Let $A^{1} =N \cap
\Gamma_{1}$ and $A^{2} =N \cap \Gamma_{2}$.
We use $a_{2} = ( 0,0,0 )$ as the base point of $A^{2}$ and $a_{1} = ( 0,0,1
)$ the base point of $A^{1}$. Let $B_{k}^{j}$ be the disk in $A^{j}$ with
boundary $C^{j}_{k}$ for $k=1,2, \ldots$ and $j=1,2$. See Figure \ref{fig:shrinking-braid}.
\begin{figure}
    \centering
    \includegraphics{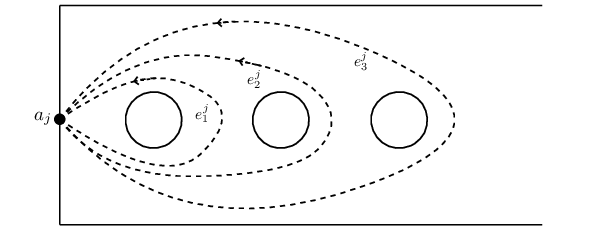}
    \caption{The dashed curves illustrate the generators of fundamental group $\pi_{1} ( \widetilde{A^{j}} )$.}
    \label{fig:generator}
\end{figure}
Denote $\widetilde{A^{j}}=A^{j} \backslash
\nobracket \bigcup_{k=1}^{\infty} B^{j}_{k}$ for $j=1,2$. For $j=1,2$ and $i\geq 1$,
we use $e_{i}^{j}$ to denote the $i$-th generator of
$\pi_{1} ( \widetilde{A^{j}} )$ drawn in Figure \ref{fig:generator}.

For $j=1,2$, the infinite free group $F_{\infty}=\langle e_{1},e_{2},\ldots \rangle$ is isomorphic to $\pi_{1} ( \widetilde{A^{j}} ,a_{j} )$ by mapping 
\begin{equation}\label{eq:identify-map}
    e_{k} \rightarrow e_{k}^{j} \; \forall k \geq 1.
\end{equation}
This gives us identification
\begin{equation}\label{eq:identify-free-group}
    \pi_{1} ( \widetilde{A^{1}} ,a_{1} )\approx F_{\infty} \approx \pi_{1} ( \widetilde{A^{2}} ,a_{2} ).
\end{equation}

Recall that $\mathcal{C}$ is the set of all cobordisms.
For each $r \in \mathcal{C}$, we denote by $N ( r )$
the exterior of $r \cup \bigcup\{B^{j}_{k}:k\geq 1, j=1,2\}$ in $N$. For $t\in [0,1]$, denote $N_{t} ( r ) =N ( r ) \cap ( \nobracket [
0, \infty ) \times [ -1,1 ] \times \{ t \}$). For $j=1,2$, denote the embedding map by $i ( r
)^{j} = \widetilde{A^{j}} \hookrightarrow N ( r )$ which induces $i ( r )^{j}_{\ast}
: \pi_{1} ( \widetilde{A^{j}} ,a_{j} ) \rightarrow   \pi_{1} ( N ( r ) ,a_{j}
)$.

\begin{theorem}\label{thm:iso-braid}
    Given $r\in \mathcal{C}$, suppose the equivalent class of $r$ is in $S$ (i.e. $r$ is a shrinking braid). Then
  $i ( r )_{\ast}^{1}$ is an isomorphism between groups $\pi_{1} (
  \widetilde{A^{1}} ,a_{1} )$ and $\pi_{1} ( N ( r ) ,a_{1} )$.
\end{theorem}
\begin{proof}
    Recall that, two shrinking braids $r_{1}$ and $r_{2}$ are said to be
    equivalent if and only if there is a diffeomorphism $G:N\rightarrow N$ such that
    $G$ is the identity map on $\partial N$ and $G ( r_{1} ) =r_{2}$. For any $w \in \mathcal{C}$, let $\widetilde{w}$ denote the equivalent class of $w$. Suppose $\widetilde{r}=\alpha(y)$ for some $y\in R$. Since $y$ can be written in the form $y=sx$ with $s \in B_{\infty}$ and $x
    \in X_{\infty}$, there is a diffeomorphism $G:N\rightarrow N$ (which is identity map on
    $\partial N$) such that $r' =G ( r )$ is the concatenation of $u',v'$ where $\widetilde{u'} =
    \alpha ( s ) $ and $\widetilde{v'} =\alpha ( x )   $. We prove theorem for $r'$
    instead of $r$. Since $G$ is a diffeomorphism preserving boundary of $N$,
    the theorem still holds for $r$. Note that $N(r')$ is the concatenation of $\overline{N ( u' )}$ and
    $\overline{N ( v' )}$ along the surface $C=N_{1/2} ( r' ) =\overline{N ( u' )} \cap \overline{N ( v' )}$. Here, $\overline{N ( u' )}$ and $\overline{N ( v' )}$ are scaling images of $N(u')$ and $N ( v' )$ along $z$ direction, respectively.
    See Figure \ref{fig:concatenation}.
    The embedding $i^{(1)}:C \hookrightarrow \overline{N ( u' )}$ induces $i_{\ast}^{( 1 )} : \pi_{1} (
    C ) \rightarrow \pi_{1} ( \overline{N ( u' )} )$ and $i^{(2)}:C \hookrightarrow \overline{N ( v' )}$ induces $i_{\ast}^{( 2 )} : \pi_{1}
    ( C) \rightarrow \pi_{1} ( \overline{N ( v' )} )$. 
    \begin{figure}
        \centering
        \includegraphics{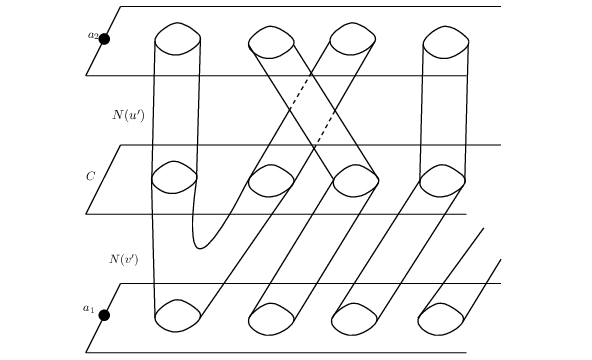}
        \caption{$r'$ is the concatenation of scaling images of $u'$ and $v'$, and $N(r')$ is the concatenation of $\overline{N(u')}$ and $\overline{N(v')}$.}
        \label{fig:concatenation}
    \end{figure}

    $i_{\ast}^{( 1 )}$ is an isomorphism by a well-known result on braids. To prove our theorem, it suffices to prove $i_{\ast}^{(2)}$ is an isomorphism. If this is true, then by using van Kampen
    theorem on $N ( r' ) =\overline{N ( u' )} \cup \overline{N ( v' )}$, we have 
    \begin{align*}
        &\pi_{1} ( N ( r' )) \approx 
        \pi_{1} ( \overline{N ( u' )})\ast_{\pi_{1} (C)} \pi_{1} ( \overline{N ( v' )}) \\
        \approx  &\pi_{1} ( \overline{N ( u' )})\approx
    \pi_{1} ( \widetilde{A^{1}})
    \end{align*}
    and our theorem follows. Here, the last isomorphism $\pi_{1} ( \overline{N ( u' )})\approx
    \pi_{1} ( \widetilde{A^{1}})$ is again by results on braids.

    Now, as in equation \eqref{eq:wordform}, we can write $x=x_{1}^{i_{1}} x_{2}^{i_{2}}\cdots x_{k}^{i_{k}}$ with integers $i_{j}\geq 0$ for $1\leq j\leq k$.

    Let $M_{l} ( l=1,2, \ldots ,k )$ be $k$ cuboids. As in definition \ref{def:j-pant}, let $\rho_{l}$ be a
    $(i_{l}+1)$-pant embedded in $M_{l}$ with exterior $M_{l}(\rho_{l}
    )$ for each $l=1,2, \ldots ,k$. Denote by $A_{l}$ the upper surface of $M_{l}$ and
    denote $\widetilde{A_{l}}=A_{l}\setminus \bigcup_{k} D^{(l)}_{k}$ where $\{D^{(l)}_{k}:1\leq k\leq i_{l}+1\}$ are input disks of $\rho_{l}$ in $A_{l}$. Now for $l=1,2, \ldots ,k-1$, we glue the right surface of $M_{l}$ and the left surface
    of $M_{l+1}$ without changing the relative position of $\rho_{l}$ inside
    $M_{l}$. Denote the resulting space $M$ and denote by
    $\rho = \coprod_{l=1}^{k} \rho_{l}$ the disjoint union of $k$ cobordisms
    embedded in $M$. Finally let $M_0$ be the resulting space by gluing the right surface of $M$ and the left surface of $N$ in which there is a trivial shrinking braid (i.e. the identity element $e$). See Figure \ref{fig:space-M_0} for an illustration. Let cobordism $\rho _0$ be the disjoint union of $\rho$ and the trivial shrinking braid (inside $M_{0}$). Denote by $M_0(\rho_0)$ the exterior of $\rho_0$ in $M_0$. Define $A_0$ to be the upper surface of $M_0$ and let $\widetilde{A_{0}}$ be the complement of input disks of $\rho_0$ in $A_0$.
    Now, $v'$ is essentially $\rho _0$. i.e. There is
    a diffeomorphism from $M_0$ to $\nobracket [
0, \infty ) \times [ -1,1 ] \times [0,\frac{1}{2}]$ such that it maps $\rho _0$ to $v'$, maps $M_{0}(\rho_{0})$ to $\overline{N(v')}$, maps $A_0$ to $\nobracket [
0, \infty ) \times [ -1,1 ] \times \{ \frac{1}{2} \}$ and maps $\widetilde{A_{0}}$ to $C$.
    \begin{figure}
        \centering
        \includegraphics{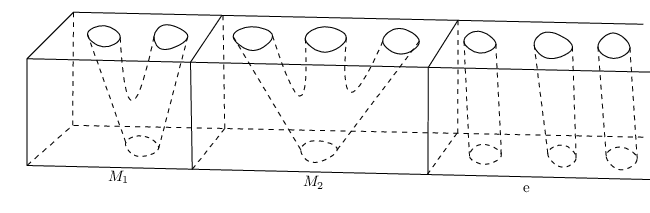}
        \caption{The picture illustrates $M_{0}$ for the case when $x=x_{1}x_{2}^{2}$.}
        \label{fig:space-M_0}
    \end{figure}

    In order to prove $i^{(2)}_{\ast}$ is an isomorphism, we only need to prove that $f'_{*}:\pi_1(\widetilde{A_{0}}) \rightarrow \pi_1(M_0(\rho_0))$ is an isomorphism where $f':\widetilde{A_{0}} \hookrightarrow M_0(\rho_0)$ is the embedding.  By
    Lemma \ref{lem:j-pant-lem}, for each $l \in \{ 1,2, \ldots ,k \}$, $i_{l} : \widetilde{A_{l}}
    \hookrightarrow M_{l} ( \rho_{l} )$ induces an isomorphism $( i_{l}
    )_{\ast} : \pi_{1} ( \widetilde{A_{l}} ) \hookrightarrow \pi_{1} ( M_{l} (
    \rho_{l} ) )$. Thus by van Kampen theorem and naturality of push out, we imply
    that $f'_{*}$ is an isomorphism and our theorem follows.
  \end{proof}
  Finally we prove Theorem \ref{thm:main}:
\begin{proof}[Proof of Theorem \ref{thm:main}]
For each $r\in S$ and $j=1,2$, recall that $i ( r )^{j} = \widetilde{A^{j}} \hookrightarrow N ( r )$
induces $i ( r )^{j}_{\ast} : \pi_{1} ( \widetilde{A^{j}} ,a_{j} ) \rightarrow
\pi_{1} ( N ( r ) ,a_{j} )$.
By \eqref{eq:identify-map} and \eqref{eq:identify-free-group}, we can regard $i ( r )^{j}_{\ast} : F_{\infty} \rightarrow
\pi_{1} ( N ( r ) ,a_{j} )$ for $j=1,2$. By Theorem \ref{thm:iso-braid}, $i ( r )^{1}_{\ast}$
is an isomorphism.
We further identify fundamental groups $\pi_{1} ( N ( r ) ,a_{1} )$ and $\pi_{1} ( N ( r ) ,a_{2} )$ which have different base points. We identify them by connecting $a_{1},a_{2}$ with the path $\gamma(t)=(0,0,1-t)$ for $t\in [0,1]$ and denote the resulting isomorphism $h:\pi_{1} ( N ( r ) ,a_{1} )\approx \pi_{1} ( N ( r ) ,a_{2} )$.
Define $\phi ( r ) = ( i ( r )^{1}_{\ast} )^{-1} \circ h^{-1} \circ (
i ( r )^{2}_{\ast} )$, then $\phi(r):F_{\infty} \rightarrow F_{\infty}$. Thus we have $\phi:S\rightarrow \tmop{Hom}(F_{\infty},F_{\infty})$.  

Since $r_{1} r_{2}$ is the
concatenation for any two shrinking braids $r_{1} ,r_{2}\in S$, we have $\phi ( r_{1} r_{2} ) = \phi ( r_{1} ) \circ
\phi ( r_{2} )$. By Proposition \ref{prop:action-R-on-F_infty} and \eqref{eq:identify-map}, we have $\phi ( \alpha ( x_{i} ) ) = \beta (
x_{i} )$ and $\phi ( \alpha ( s_{i} ) ) = \beta ( s_{i} )$ for each $i \in
\mathbbm{N}$. Since $R$ is generated by $x_{i},s_{i}$'s, we have $\phi \circ \alpha = \beta$. By Theorem \ref{thm:beta-injective}, $\beta$ is injective and thus $\alpha$ is injective. Together with Proposition \ref{prop:define-alpha}, $\alpha$ is bijective and Theorem \ref{thm:main} follows.
\end{proof}













\bibliographystyle{hal}
\bibliography{bib}

\end{document}